\documentclass{amsart}
          
\usepackage{fullpage,eucal,amsmath,amsthm,amsfonts,verbatim,amsopn,amssymb}

\usepackage[all]{xy}

\newcommand\dmo{\DeclareMathOperator}

\renewcommand\sec{\section}
\newcommand\subs{\subsection*}
\newcommand\para{\paragraph*}
\newcommand\bea{\begin{array}{r@{,\ \ \ \ \ \ \ }l}}
\newcommand\ena{\end{array}}
\newcommand\barr{\begin{array}}
\newcommand\earr{\end{array}}
\newcommand\ben{\begin{enumerate}}
\newcommand\een{\end{enumerate}}
\newcommand\beq{\begin{equation}}
\newcommand\eeq{\end{equation}}
\newcommand\bqa{\begin{eqnarray*}}
\newcommand\eqa{\end{eqnarray*}}
\newcommand\bqan{\begin{eqnarray}}
\newcommand\eqan{\end{eqnarray}}
\newcommand\bit{\begin{itemize}}
\newcommand\eit{\end{itemize}}
\newcommand\bfi{\begin{figure}[htbp]}
\newcommand\efi{\end{figure}}
\newcommand\bce{\begin{center}}
\newcommand\ece{\end{center}}

\newcommand\bpr{\begin{proof}}
\newcommand\epr{\end{proof}}

\newcommand{\ignore}[1]{}

\newtheorem{theorem}{Theorem}[section]
\newtheorem{conjecture}[theorem]{Conjecture}
\newtheorem{lemma}[theorem]{Lemma}
\newtheorem{corollary}[theorem]{Corollary}
\newtheorem{proposition}[theorem]{Proposition}

\theoremstyle{definition}

\newtheorem{definition}[theorem]{Definition}
\newtheorem{example}[theorem]{Example}
\newtheorem{remark}[theorem]{Remark}
\newtheorem{remarks}[theorem]{Remarks}

\newtheorem{question}[theorem]{Question}

\newcommand\bthm{\begin{theorem}}
\newcommand\ethm{\end{theorem}}
\newcommand\bcn{\begin{conjecture}}
\newcommand\ecn{\end{conjecture}}
\newcommand\bla{\begin{lemma}}
\newcommand\ela{\end{lemma}}
\newcommand\bco{\begin{corollary}}
\newcommand\eco{\end{corollary}}
\newcommand\bpro{\begin{proposition}}
\newcommand\epro{\end{proposition}}
\newcommand\bdf{\begin{definition}}
\newcommand\edf{\end{definition}}
\newcommand\bex{\begin{example}}
\newcommand\eex{\end{example}}
\newcommand\brm{\begin{remark}}
\newcommand\erm{\end{remark}}
\newcommand\brms{\begin{remarks}}
\newcommand\erms{\end{remarks}}
\newcommand\bqu{\begin{question}}
\newcommand\equ{\end{question}}

\newcommand\msk{\medskip}

\newcommand\prf{\textsc{Proof. \quad}}

\newcommand\sr{\stackrel}
\newcommand\<{\langle}
\renewcommand\>{\rangle}
\newcommand\ot{\otimes}
\newcommand\op{\oplus}

\newcommand\sub{\subseteq}

\newcommand\ra{\rightarrow}

\newcommand\lra{\longrightarrow}

\newcommand\seq{\sim}
\newcommand\teq{\approx}
\newcommand\teqq{\parallel}

\newcommand\al{{\alpha}}
\newcommand\be{{\beta}}

\newcommand\dl{\delta}
\newcommand\dltl{\tilde{\delta}}
\newcommand\lm{\lambda}
\newcommand\Lm{\Lambda}

\newcommand\om{\omega}
\newcommand\omt{\omega^\tau}
\newcommand\Om{\Omega}
\newcommand\Omt{\Omega^\tau}

\newcommand\pd{\partial}

\newcommand\xbar{\overline{x}}

\newcommand\abar{\overline{a}}
\newcommand\bbar{\overline{b}}

\newcommand\A{{\mathbb{A}}}

\newcommand\Z{{\mathbb{Z}}}

\renewcommand\P{{\mathbb{P}}}

\renewcommand\O{{\mathcal{O}}}

\newcommand\aff{\textup{Aff}}
\newcommand\Hom{\textup{Hom}}

\newcommand\Ext{\textup{Ext}}
\renewcommand\dim{\textup{dim}}
\newcommand\HH{\textup{H}}
\newcommand\hh{h}
\newcommand\Spec{\textup{Spec}}

\newcommand\Ker{\textup{Ker}}

\newcommand\acl{\textup{acl}}

\dmo{\spec}{Spec}

\newcommand\sm{strongly minimal\,}

\newcommand\prl{^{(1)}}

\newcommand{\GL}{\textup{GL}}

\newcommand{\Aff}{\textup{Aff}}
\newcommand{\Tr}{\textup{Tr}}

\newcommand\ME{{\mathcal{E}}}
\newcommand\MF{{\mathcal{F}}}
\newcommand\MG{{\mathcal{G}}}
\newcommand\MH{{\mathcal{H}}}
\newcommand\MFT{{\mathcal{F}^\tau}}
\newcommand\MGT{{\mathcal{G}^\tau}}
\newcommand\tis{$\tau$-invertible sheaf}
\newcommand\tiss{$\tau$-invertible sheaves}
\newcommand\tauf{$\tau$-form}
\newcommand\taufs{$\tau$-forms}
\newcommand\PO{PO}

\begin{document}

\title{Order 1 strongly minimal sets in differentially closed fields}
\author{Eric Rosen}
\address{Department of Mathematics\\
Massachusetts Institute of Technology\\
Cambridge, MA 02139}
\email{rosen@math.mit.edu}
\urladdr{http://math.mit.edu/~rosen/}
\thanks{This paper constitutes a revised version of part of my dissertation, at the 
University of Illinois at Chicago (2005).  I am grateful to Lawrence Ein,
Henri Gillet, David Marker, and Thomas Scanlon for helpful discussions on this material.}

\date{\today}

\begin{abstract}
We give a classification of non-orthogonality classes of trivial order 
1 strongly minimal sets in differentially closed fields.  A central idea 
is the introduction of $\tau$-forms, functions on the prolongation of 
a variety which are analogous to 1-forms.  Order 1 strongly minimal sets then 
correspond to smooth projective curves with $\tau$-forms.  We also
formulate our results scheme-theoretically, in terms of $\tau$-differentials
and $\tau$-invertible sheaves on curves, thereby obtaining additional 
information about the strongly minimal sets.
This work partially generalizes and extends results of Hrushovski and Itai. 
\end{abstract}

\maketitle

\sec*{Introduction}

This paper addresses the problem of classifying the strongly minimal sets
definable in differentially closed fields.  A nice description of 
the background to the question and known results can be found in the 
introduction of \cite{HI},
so we only give a brief introduction here.  For more on the model theory of 
differential fields, see also \cite{MMP}, \cite{Pil02}, 
\cite{HS}, and \cite{HICM}.
The precise problem is to classify strongly minimal\ sets up to non-orthogonality, 
a natural logical equivalence relation.  Associated to each strongly minimal\ set
is a combinatorial geometry, which is either non-locally modular, locally
modular non-trivial, or trivial, thus dividing strongly minimal\ sets into three
classes.  Hrushovski and Sokolovi\'c \cite{HS} classified
both the non-locally modular and the locally modular non-trivial strongly minimal\ sets.
The former are non-orthogonal to the field of constants, which is strongly minimal,
and the latter correspond to (isogeny classes of) simple abelian varieties 
that do not descend to the constants.  What remains is to classify the 
trivial ones.

Let $K$ be a universal domain for DCF, the theory of differentially closed
fields, and let $k$ denote the field of constants.
A strongly minimal\ set, or formula, determines a unique strongly minimal\ type, and sets/formulas 
determining the same type are non-orthogonal, so one can equally well consider 
strongly minimal\ types.  Recall that if $q$ is a strongly minimal\ type over a differential field $L$, 
then the {\em order} of $q$ is $\textup{td}((\abar, \abar', \ldots)/L)$,
for $\abar$ realizing $q$, which is always finite.  
It seems to be unknown whether there are (trivial)
strongly minimal\ types of arbitrarily large order, though there are examples of order 2.
One consequence of \cite{HI} is that there are many types of order 1.

Hrushovski \cite{HJ} proved that every order 1 strongly minimal\ type that is orthogonal to $k$
is $\om$-{\em categorical} and thus trivial, since the classification of 
the non-trivial locally modular strongly minimal\ types by Manin kernels implies that
they are never $\om$-categorical.  
(Alternatively, one can derive the triviality of 
order 1 locally modular strongly minimal sets 
as a consequence of the fact that any non-trivial locally modular 
strongly minimal set is non-orthogonal to the Manin kernel of a simple abelian
variety~\cite{HS} and Buium's proof that any such Manin kernel has transcendence 
degree at least 2).  With Itai, Hrushovski then gives a quite precise
description of the trivial strongly minimal\ types over $k$, which we summarize below.
In this paper, we consider all order 1 strongly minimal\ types,
partially generalizing and extending some results of \cite{HI}.

\subs{Background}
It is well-known, and easy to see, that (complete, stationary) $n$-types in
pure algebraically closed fields correspond to irreducible Zariski closed 
subsets of $\A^n$.  Call such a set an {\em embedded affine variety},
as it is an affine variety together with a closed embedding into affine space.
We recall a similar  geometric description of complete strongly minimal\ types in DCF, 
which also holds more generally for types of finite rank.  
(See \cite{HI} or \cite{Pil02} for more details.)  

Say that two types $p$ and $q$ are {\em interdefinable}
if, perhaps passing to non-forking extensions over a common base field,
there is a definable bijection between their sets of realizations.  Then
any strongly minimal\ type $p$ in DCF is interdefinable with a type $q$ such that for
any realization $\abar$ of $q$, $\dl(\abar) = s(\abar)$, where $s(\xbar)$ is a
tuple of polynomials.  Since interdefinable types are non-orthogonal, we
may assume that all (strongly minimal) types are of this form.  

Given an affine variety $V \sub \A^n$, the map $\dl$ above, acting 
componentwise, is a section of the first prolongation $\tau V \sub
\A^{2n}$, a (possibly reducible) affine variety defined below,
also called the {\em shifted tangent bundle}.
By the above reduction, a strongly minimal\ type $p$ corresponds to a pair
$(V,s)$, where $V \sub \A^n$ is an embedded affine variety,
and $s: V \ra \tau V\sub \A^{2n}$ 
is a section of the {\em first prolongation}.  In this case,
one says that $p$ {\em lives on} the variety $V$.  It is easy to check
that the order of $p$ equals the dimension of $V$, so 
the order 1 types are exactly those that live on curves.

In general, a pair $(V,s)$ as above determines a Kolchin closed set
of finite Morley rank, $\Xi(V,s):= \{v \in V | \dl(v)=s(v) \}$.  When $V$
is a curve, this set is necessarily strongly minimal, so there's actually a bijection
between (Kolchin closed) order 1 strongly minimal sets and pairs $(C,s)$, $C$ an embedded
affine curve and $s$ a section of the prolongation, $s:C\ra \tau C$.
(Moreover, this bijection preserves fields of definition.)
This gives a precise, purely geometric characterization of order 1
strongly minimal sets.
What remains, though, is to determine whether or not such a 
strongly minimal\ set $\Xi(C,s)$ is trivial, and to understand non-orthogonality.

If $X$ is a strongly minimal\ set living on a curve $C$, and $C$ is birational to
some $C'$, then $X$ is interdefinable with some $X'$ living on $C'$.
Since every curve is birational to a unique smooth projective curve,
it suffices to consider such sets living on (embedded) smooth projective
curves.  Strictly speaking, then, we want to classify strongly minimal
sets of the form $\Xi(C,s)$, $C$ an embedded smooth projective
curve, $s: C \ra \tau C$.  But we will consider $C$ as an
abstract curve, since the particular embedding in $\P^n$ is unimportant,
as different embeddings yield interdefinable strongly minimal sets.

Given a curve $C$ (or any variety) defined over $k$, the prolongation 
$\tau C$ equals the tangent variety $TC$.  In this case, a section
of the prolongation $s: C \ra \tau C$ is just a vector field. For curves, 
there is also a natural bijection between vector fields and 1-forms.
Given a curve $C$ and a vector field $s$, let $\om$ be the 1-form such
that $\om (s) = 1$ almost everywhere.  Thus, order 1 strongly minimal\ sets defined over
the constants are also represented by pairs $(C,\om)$, defined over
the constants.  

Using this representation, Hrushovski and Itai obtain a rather complete 
description of  non-orthogonality classes of trivial order 1 strongly minimal\ sets 
defined over the constants.  To any such class they associate a unique 
(smooth projective) curve, the main idea being to pick a canonical 
{\em strictly minimal} set in any such class.  Second, they prove that for 
any curve of genus $\geq 2$, there are many classes associated to that curve, 
thereby proving that there are indeed `very many' trivial order 1 strongly minimal\ sets 
defined over the constants.

\subs{Results}

The original motivation for this work was to generalize results from
\cite{HI} to all order 1 strongly minimal\ sets.  One of the main new ideas is the 
introduction of $\tau$-forms, which are functions on the prolongation of
a curve, analogous to 1-forms, which are functions on the tangent variety.
This enables one to develop a `geometric approach' to order one strongly minimal\ sets
as in \cite{HI}, with a pair $(C,\omt)$, $C$ a curve, $\omt$ a $\tau$-form, 
representing a strongly minimal\ set.  We can then characterize non-orthogonality
classes of trivial strongly minimal\ sets in terms of `essential $\tau$-forms', 
appropriately defined, generalizing a central result of \cite{HI}.

In a separate paper~\cite{Ros}, we introduce and develop the theory of $\tau$-differentials,
the algebraic analog of $\tau$-forms. We use this material here to give an alternative
characterization of order one strongly minimal sets in terms of curves with
$\tau$-invertible sheaves.  
This makes it possible to determine the dimension of the space of global 
$\tau$-forms on a curve, in terms of its Kodaira-Spencer class.  
Further, we show that the set of global $\tau$-forms on a curve is (uniformly) definable,
as is the dimension of this set.
It also leads to the introduction of the prolongation cone of a curve, a vector bundle
into which both the tangent variety and prolongation naturally embed.

\subs{Open questions}

The main question left open by this work is whether for every curve $C$
of genus $>1$ there are $\sim$-essential \taufs.  For curves defined over
$k$, a positive answer follows from \cite{HI}.  A positive
answer for all curves would show that there are order 1 trivial strongly minimal\ sets 
that are orthogonal to any such set defined over $k$.  

Two main problems about trivial strongly minimal\ sets in DCFs, asked by
Hrushovski, remain completely open.  Are they all $\om$-categorical?
And can they be classified?  Extending the classification beyond order 1
certainly seems difficult.  For example, presumably one would need to 
be able to determine whether a finite rank type given as a pair $(V,s)$ 
is strongly minimal.  But perhaps there are particular cases of order 2 types,
thus living on surfaces, that are more accessible.  
There are also related questions, which might be approached
by existing methods.  Hrushovski suggests that one could try to show
that the solutions to a generic order 2 equation are strongly minimal, 
to answer a question of Poizat.  
Another well-known question is whether every U-Rank 1 type is strongly minimal.  
This is related to Hrushovski and Scanlon's result~\cite{HrSc} 
that U-Rank and Morley rank differ in DCF.

\sec{Prolongations}
\label{prolong}
We recall the construction of the prolongation of a variety over a differentially
closed field.  Buium's original definition uses the language of schemes, 
though we give a description in local coordinates used more often by
model-theorists.  For a variety $V$ defined over the constants, its prolongation
is just the tangent variety $TV$.  In general, though, it is a $TV$-{\em torsor}.
Thus it is also called the {\em shifted tangent bundle}.  
In this case, the fibers of the canonical projection to $V$ are no longer
vector spaces, but affine spaces, defined below.

\para{Affine Spaces}  
An affine space is essentially a vector space without a 
distinguished point as origin.  Affine spaces arise in model theory as combinatorial
geometries that are locally modular but not modular.   Just as a locally
modular geometry becomes modular when one fixes any point, choosing a point in 
an affine space naturally gives a vector space.  The presentation here
differs from than that in, for example, Hodges' textbook~\cite{Hod}, 
but is basically 
equivalent.  (The word ``affine'' will also be used in a 
completely different sense, in connection with affine varieties.)

\bdf
\label{aff}
Let $K$ be a field.  A {\em K-affine space} is a triple $(A, V, \alpha)$,
where $A$ is a set, $V$ is a $K$-vector space, and $\alpha$ is a regular
action of $V$ on $A$. The {\em dimension} of $A$, $\dim(A)$, is just $\dim(V)$.
\edf

We will write $\alpha(v,a)$ as $v \cdot a$, and often omit the action $\alpha$ 
when it is understood.

\brm
\label{affbasics}
\ben
\item  For any $K$-vector space $V$, $(V, V)$ is an affine space in a natural way.
\item  Given $(A, V)$ and $a \in A$, there is a natural bijection $i_a:V\lra A$
given by $i_a(v) = v\cdot a$.
\een
\erm

\bdf
An {\em affine map} between $K$-affine spaces $(A, V, \alpha)$ and $(B, W, \beta)$
is a map $f:A \lra B$ such that there is a linear map $\lambda f: V \lra W$ so that 
for all $a\in A, v \in V$, $f(v\cdot a) = (\lambda f (v))\cdot f(a)$.
Let $\Aff(A,B)$ denote the set of all affine maps from $A$ to $B$.
(By putting coordinates on $A$ and $B$, this set can be endowed with a vector
space structure.)

\edf
   
\brm
\label{affrems}
\ben
\item  Given an affine map $f$, the linear map $\lambda f$ is uniquely determined.
In fact, $\lambda$ is a functor from the category of $K$-affine spaces to the 
category of $K$-vector spaces.
\item Given affine maps $f, g$, from $(A, V)$ to $(B, W)$, $\lambda f = \lambda g$ 
if and only if there is a $w \in W$ such that for all $a \in A$, $f(a) = w \cdot g(a)$.
\item There is a short exact sequence 
$$
0 \lra B \stackrel{\mu}{\lra}\aff(A, B)\stackrel{\lambda}{\lra}\Hom(V, W) \lra 0
$$
where for all $b_0 \in B, \mu(b_0)$ is the constant affine map 
$\mu(b_0): A \lra B$ mapping each $a \in A$ to $b_0$.  Thus $\dim(\aff(A, B)) = 
\dim(B) + \dim(A)\dim(B)$.
\een
\erm

We are particularly interested in the case, $W = K$, i.e., when $B$ is 
1-dimensional.
Then $\aff(A, B)$ is something like a dual space to $(A, V)$, though 
$\dim(\aff(A,B))=\dim(A) + 1$.    

\brm
Let $(A,V)$ be an affine space.  Recall that the {\em affine group} of $V$,
denoted $\Aff(V)$, is the group generated by the translation group $\Tr(V)$
and $\GL(V)$.  In fact, $\Aff(V)=\Tr(V)\rtimes \GL(V)$.
The automorphism group of $(A,V)$ is naturally isomorphic to $\Aff(V)$.
\erm

\para{Varieties, tangent spaces, and prolongations}
A prevariety is a topological space with an open cover,
$W = W_1 \cup \ldots \cup W_m$ and compatible coordinate charts, $f_i:W_i \lra V_i$,
$V_i$ an affine variety.  For all $i, j  \leq m$, let $U_{i, j} = f_i(W_i \cap 
W_j) \sub V_i$ and $f_{i, j}: U_{i, j} \lra U_{j, i}$ be $f_j\circ f_i^{-1}$.
For our purposes, a variety will be an irreducible, smooth separated prevariety.

Recall that, given a polynomial $p(X)$, 
$p^\delta(X)$ denotes the polynomial obtained by applying $\delta$ to
each of the coefficients.  For more information, see 
Marker~\cite{Mar03} or Pillay~\cite{Pil02}.  The following easy observation 
will be useful.

\bla
For any $n$, let $\epsilon: K[x_1, \ldots , x_n] \lra K[x_1, \ldots , x_n]$
be the map that takes any polynomial $f$ to $f^\dl$, which is obtained 
from $f$ by taking the derivative of each coefficient.  Then $\epsilon$
is a derivation on $K[x_1, \ldots , x_n]$.  Further, $\epsilon$ commutes
with each derivation $\frac{d}{dx_i}$. 
\ela

\bpr
Straightforward.
\epr

\bdf
Let $V \subseteq K^n$ be an irreducible smooth affine variety, 
$I(V) = \<p_1,\ldots , p_m\>$.  The tangent space of $V$ is
$$
TV = \{(a, u) \in K^{2n} : \sum_{i=1}^n\frac{\partial p_j}{\partial x_i}(a) 
\cdot
u_i = 0, j = 1, \ldots , m\}
$$
The first prolongation of $V$ is 
$$
\tau V = \{(a, u) \in K^{2n} : \sum_{i=1}^n\frac{\partial p_j}{\partial x_i}(a) 
\cdot u_i  + p_j^\delta(a)= 0, j = 1, \ldots , m\}
$$

Since everything is functorial, one can also define $TV$ and $\tau V$
for general varieties, in a coordinate free manner.
\edf

\brm
There are natural projection maps $\pi_T: TV \lra V$ and $\pi_\tau: \tau V \lra
V$.  For smooth $V$, and $v \in V$, $\pi_T^{-1}(v)=TV_a$ is a $\dim(V)$-dimensional
vector space, and $\pi_\tau^{-1}(v)=\tau (V)_a$ is a $\dim(V)$-dimensional 
affine space naturally acted on by $\pi_T^{-1}(v)$.

Further, let $TV \times_V\tau V = \{(a, u, w) \in K^{3n} : (a, u) \in TV, (a, w) \in 
\tau V\}$.  The map $p: TV \times_V\tau V \lra \tau V$ given by
$p(a, u, w) = p(a, u+w)$ makes $\tau V$ a {\em torsor} under $TV$.  
\erm

We now introduce a new construction, the prolongation cone of a variety,
into which both the tangent variety and prolongation naturally embed.  For more
details, see~\cite{Ros}.

\bdf
Let $V \subseteq K^n$ be an irreducible smooth affine variety, 
$I(V) = \<p_1,\ldots , p_m\>$.  The {\em prolongation cone} of $V$ is
$$
\tilde{V} = \{(a, u, u') \in K^{2n+1} : \sum_{i=1}^n\frac{\partial p_j}{\partial x_i}(a) 
\cdot
u_i + p_j^\delta(a)\cdot u'= 0, j = 1, \ldots , m\}
$$

As in the case of the tangent variety, there is a  natural projection map
$\tilde{\pi}:\tilde{V} \lra V$, whose fibres are vector spaces, now of dimension
$\dim (V) + 1$.  More generally, the prolongation cone is a $\dim(V) + 1$ vector
bundle over $V$.  Again, the construction also globalizes to general varieties.

Observe that the intersection of the prolongation cone with the hyperplane $u' = 0$
is the tangent variety.  Likewise, the intersection with the hyperplane $u' = 1$
is the prolongation.
\edf

\brm  
The tangent space $TV$ is a variety with `additional linear structure' on the fibers 
$TV_v$.  We recall
how to make this notion precise and indicate how to make rigorous the 
notion that the fibers $\tau V_v$ of the prolongation are affine spaces.

One can define a vector bundle as follows (see \cite{GHL}, p.\ 15).
Let $E$ and $V$ be smooth varieties over a field $K$, $\pi:E\lra V$ a regular map.
We say that $(\pi , E, V)$ is a {\em vector bundle of rank n} if the following
conditions hold.
\ben
\item  $\pi$ is surjective.
\item  There exists a finite open cover $(U_i)_{i\in I}$ of $V$, and
  isomorphisms $h_i:\pi^{-1}(U_i)\lra U_i\times K^n$ such that for any
  $x\in U_i$, $h_i(\pi^{-1}(x))=\{x\}\times K^n$.
\item  For any $i,j \in I$, there is a regular map $g_{ij}:U_i\cap U_j \lra
\GL_n(K)$ such that the map
$$
h_i \circ h_j^{-1}:(U_i\cap U_j) \times K^n \lra (U_i\cap U_j) \times K^n 
$$
is of the form $h_i\circ h_j^{-1}(x,v)=(x,g_{ij}(x)\cdot v)$.
\een
  
$V\times K^n$ is the {\em trivial vector bundle} over $V$, and the definition 
of a vector bundle says that it is {\em locally trivial}.  Likewise, the 
prolongation $\tau V$ is a locally trivial affine bundle, in a similar sense.  
That is, let $V$
be a variety and $(K^n,K^n)$ $n$-dimensional affine space.  Then 
$(\pi , E, V)$ as above is a {\em trivial affine bundle of rank n} if it is 
isomorphic to $V\times K^n$, where we consider $K^n$ as an affine, rather than
vector, space.  More generally, $(\pi , E, V)$ as above is
an {\em affine bundle of rank n} if we have the following.
\ben
\item  $\pi$ is surjective.
\item  There exists a finite open cover $(U_i)_{i\in I}$ of $V$, and
  isomorphisms $h_i:\pi^{-1}(U_i)\lra U_i\times K^n$ such that for any
  $x\in U_i$, $h_i(\pi^{-1}(x))=\{x\}\times K^n$.
\item  For any $i,j \in I$, there is a map $g_{ij}:U_i\cap U_j \lra
\Aff(K^n)$ such that the map
$$
h_i \circ h_j^{-1}:(U_i\cap U_j) \times K^n \lra (U_i\cap U_j) \times K^n 
$$
is of the form $h_i\circ h_j^{-1}(x,v)=(x,g_{ij}(x)\cdot v)$.
\een

\erm

\brm 
Another way to formalize the notion of an affine bundle is suggested by
the notion of a fiber bundle associated to a principal bundle
(See \cite{KN96}, p.\ 50--55.  Compare also section III.3 on 
Affine Connections, p.\ 125.)  
Roughly, given a manifold $M$ and a Lie group $G$,
a principal $G$-bundle $P(G,M)$ is a fiber bundle 
$\pi:E\lra M$, with a free $G$-action
on each fiber $\pi^{-1}(x), x\in M$.  If $F$ is some other manifold, with 
a $G$-action, then one can construct a fiber bundle $E(G,M,F,P)$ which is
the fiber bundle over $M$, with standard fiber $F$, and structure group
$G$, associated to the principal fiber bundle $P$.  In our setting, $M$ is our 
variety $V$, $G$ is $\Aff(K^n)$, and $F$ is our affine space.
\erm

The following lemma is due to Buium \cite{Bui93}.
\bla
Let $V$ be a smooth variety of dimension $n$.  Then $\tau V$ is an affine bundle
of rank $n$.
\ela

\para{Tangent and lifting maps}
We now consider maps between varieties.  Let $V,W$ be varieties,
and $\phi:V \ra W$ a regular map.  The map $\phi$ determines a map
from $TV$ to $TW$, written $T\phi$,  the {\em tangent map} of $\phi$,
which restricts, for each $a \in V$,
to a linear map on fibers, $T\phi_a:TV_a \ra TW_{f(a)}$.

An important special case occurs when $f$ is a regular function on $V$,
viewed as a map from $V$ to $\A^1$.  We consider the {\em differential},
$df$, which is the composition $df = \pi\circ Tf$, where 
$\pi$ is the projection from $T\A^1 \ra \A^1$ onto the {\em tangent vector
component}.  Thus, $df$ is a regular function on $TV$, and we have a
$K$-linear derivation $d:K[V]\ra K[TV]$.  Alternatively, $d$ is a
derivation from $K[V]$ to $\Om[V]$, the regular differential forms
on $V$, which is a $K[V]$-module with a natural embedding in $K[TV]$.

For prolongations, there is also a 
{\em lifting map} from $\tau V$ to $\tau W$, 
which we write $\phi^{(1)}$, which restricts, for each $a\in V$, to an 
affine map on fibers, $\phi^{(1)}_a:\tau V_a \ra \tau W_{f(a)}$.
For affine varieties, $V \sub K^n, W \sub K^m$, the map $\phi^{(1)}$
is given by $\phi^{(1)}_a((a,u)) = (\phi(a),d\phi_a(u)+\phi^\dl(a))$.

When $f$ is a regular function on $V$, composing $f^{(1)}$ with 
$\pi$, as above, one gets a map $\tau f = \pi\circ f^{(1)}$,
which we call a {\em $\tau$-differential}.
Below, Lemma~\ref{twder},
we will see that $\tau:K[V]\ra K[\tau V]$ is a derivation (in fact,
a $\tau$-derivation, as defined below.

\brms
\ben
\item  When the derivation $\dl$ on $K$ is trivial, then for any variety
$V$, $\tau V = TV$.  Likewise, given any regular map between varieties,
$f:V \ra W$, $f\prl= Tf$.  More generally, this holds true over an
arbitrary differential field if everything is defined over the field of 
constants.

\item  In \cite{Ros}, Proposition 3.15, it is shown that that the map 
$\tau$ coincides with something introduced by Buium, 
in a different context. 
\een
\erms

The next lemma follows from the description, above, of the 
lifting map in local coordinates.

\bla
\label{taudiff}
Let $V$ be a variety, $f$ a regular function on $V$.  For each 
$a\in V$, $\lm(\tau f_a)=df_a$.
\ela


\sec{$\tau$-forms}
\label{tauforms}  
Our main idea is to look at an analog of 1-forms on $\tau V$, which will be functions
$f:\tau V \lra K$ such that on each fiber of $\pi_\tau : \tau V \lra V$, $f$ is 
an affine map, as defined above.
We recall the definition of 1-forms, following the presentation of 
Shafarevich~\cite{Sha}.  We then use this formalism
to introduce $\tau$-forms.   
For an alternative, more general treatment of this material, 
in the language of commutative algebra and scheme theory, see \cite{Ros}.

Let $V$ be a variety, and let $\Phi [V]$ be the set of all functions 
$\phi$ mapping each point $v \in V$ to a linear map 
$\phi (v): TV_v\lra K$.
Note that $\Phi [V]$ naturally forms a (very large) $K[V]$-module. 
Given $f \in K[V]$, the differential $df$ is a function in $\Phi [V]$.
One could look at the submodule of $\Phi[V]$ generated by
$\{df: f \in K[V]\}$, but this is somewhat too small.
Instead, we say that an element $\phi \in \Phi [V]$ is a regular differential
1-form on $V$ if every $v\in V$ has a neighborhood $U$ such that the 
restriction of $\phi$ to $U$ belongs to the $K[U]$-submodule of 
$\Phi[U]$ generated by the elements $df, f \in K[U]$.

The regular differential 1-forms on $V$ form a $K[V]$-module, 
denoted $\Om [V]$.

\bla
The map $d: K[V]\lra\Om[V]$ satisfies 
\ben
\item  $dc=0$ for all $c \in K$.
\item  $d(f + g) = df+dg$
\item  $d(fg)= fdg + gdf$.
\een
\ela

\bpro
Every $v\in V$ has an affine neighborhood $U$ such that 
$\Om [U]$ is a free $K[U]$-module of rank $\dim (V)$.
\epro

We now introduce $\tau$-forms, imitating this construction.
Let $\Psi[V]$ be the set of all functions $\psi$ mapping each 
$v \in V$ to an affine map $\psi(v): \tau V_v\lra K$.
As above, $\Psi [V]$ forms a $K[V]$-module and,
given $f \in K[V]$, the $\tau$-differential $\tau f$ is in $\Psi [V]$.  
  
\bdf
We say that $\psi\in \Psi[V]$ is a {\em regular $\tau$-form} if every $v\in V$ 
has a neighborhood $U$ such that the restriction of $\psi$ to $U$ belongs 
to the $K[U]$-submodule of $\Psi[U]$ generated by the elements
$\tau f,f \in K[U]$.  

The regular $\tau$-forms on $V$ form a $K[V]$-module, 
denoted $\Omt [V]$.
\edf

Note that for any variety $V$, \tauf\ $\omt \in \Omt [V]$,
and $v \in V$, there is an open neighborhood $U$ of $v$ such 
that locally, on $U$,
$$
\omt = (\sum_{i=1}^ng_i\tau f_i)
$$
for $g_i,f_i\in K[U]$.

\brm
Given a variety $V$, note that a differential form on $V$ is a 
regular function on $TV$.  Thus the module $\Om[V]$ embeds
naturally in $K[TV]$, and the differential map is a $K$-linear derivation,
$d: K[V]\ra K[TV]$.  Likewise, a $\tau$-form on $V$ is a regular
function on $\tau V$.  By Lemma~\ref{twder}, below, $\tau$ is actually
a derivation.
\erm

\brm
\label{iotaiota}
From the definition, there is a map $\tau: K[V]\lra\Omt[V]$, much
like the derivation map $d: K[V]\lra\Om[V]$ to 1-forms.  
But for \taufs\ there is also a natural embedding of $K[V]$ into
$\Omt[V]$, which we now describe.
Because $K$ is differentially closed, there is a $c\in K$ 
such that $\dl(c)=1$.  So $c\in K[V]$ and $\tau c\in \Omt[V]$ is
the constant function on $\tau V$ with value 1.  By the definition
of $\Omt[V]$, for each $f\in K[V], f\tau c \in \Omt[V]$, where
$f\tau c$ is constant on each fiber $\tau V_v$, with value $f(v)$.
One sees immediately that the map $\iota:K[V]\lra \Omt [V]$, given by
$\iota(f)=f\tau c$, is an embedding (and is independent of the choice
of $c$).  Occasionally, given $f \in K[V]$, we will also write $f \in \Omt[V]$, 
where to be more precise we mean that $\iota f\in  \Omt[V]$.

Call a \tauf\ $\omt$ {\em trivial} if it equals $\iota f$, for some $f\in K[V]$.
(Clearly, a \tauf\ is trivial if and only if it is trivial
on some non-empty open subset.)  Further, given $e\in K \sub K[V]$,
call $\iota e$ a {\em constant trivial} $\tau$-form.
\erm

\bpro
\label{freemod}
Every $v\in V$ has an (affine) neighborhood $U$ such that 
$\Omt [U]$ is a free $K[U]$-module of rank $\dim (V) + 1$.
\epro

\prf 
%
We know that $v\in V$ has a neighborhood $U_1$ such that 
$\Om[U_1]$ is a free $K[U_1]$-module of rank $\dim (V)$.
Clearly, there is a neighborhood $U_2$ of $v$ such that there is a 
regular section $E_2$ of $\tau U_2$.  In other words, there is 
a map $\phi_2: U_2 \lra \tau U_2$ that is an isomorphism from 
$U_2$ to $E_2$, whose inverse is the projection map 
$\pi_\tau|_{U_2}$.  Let $U=U_1\cap U_2$, $E=E_2\cap\tau U$,
and $\phi =\phi_2|_U$, so that $E$ is a section of $U$.

$\tau U$ is a $TU$-torsor, and the section $E$ gives us an
isomorphism $\phi'$ from $TU$ to $\tau U$, which maps the 
0-section of $TU$ to $E$.  Given $u\in TU_a$,
$$
\phi'(u)=u\cdot\phi(\pi_T(u))
$$
where $u\cdot\phi(\pi_T(u))$ denotes the action of $TU_a$
on $\tau U_a$.
(On each affine fiber, $(\tau V_a, TV_a)$, this is the 
bijection from $TV_a$ to $\tau V_a$ that one gets by
fixing the point $\phi(a) \in \tau V$, as in Remark~\ref{affbasics}.2.)

The isomorphism $\phi'$ determines a bijection $\Phi'$ between 
$\Om[U]$ and the set $\Psi\sub\Omt[U]$ of \taufs\ on $U$ 
that take the value 0 on all of $E$.  Given $\om\in\Om[U]$, let
$\Phi'(\om)$ be the \tauf\ on $U$ such that for $u \in \tau U$,
$$
\Phi'(\om)(u) = \om(\phi'^{-1}(u)).
$$
In fact, this bijection is an isomorphism of $K[U]$-modules.

We want to show that $\Omt[U]=\Psi \oplus \iota K[U]$.
Let $\omt\in\Omt[U]$, and let $g\in K[U]$ be the function
$g(a)=\omt(\phi(a))$.  Then $\omt_0=\omt -\iota g$ is in 
$\Psi$, which implies that $\Omt[U]=\Psi + \iota K[U]$.

Finally, observe that $\Psi + \iota K[U]$ is a free direct
sum $\Psi \oplus \iota K[U]$, making $\Omt[U]$ into a free $K[U]$-module of dimension
$\dim(V)+1$.  For suppose that there is an $\omt \in \Psi,
\omt=\Phi'(\om)$, $\om\in\Om[U]$, and a $g\in K[U]$ such that 
$\omt=\iota g$.  Then $\omt$ is constant on each affine $\tau$-fiber,
so $\om$ must be the trivial 1-form, and $\omt$ must be identically
zero.  Thus, $g$ is also identically zero, as desired.
\qed

\msk

The following lemma is suggestive.

\bla
\label{twder}
The map $\tau: K[V]\lra\Omt[V]$ satisfies 
\ben
\item  $\tau c=\dl c$ for all $c \in K$.
\item  $\tau(f + g) = \tau f + \tau g$
\item  $\tau(fg)= f\tau g + g\tau f$.
\een
\ela

\prf  Conditions 1.\ and 2.\ are immediate from the definitions.
To prove 3., we show first that 
$(fg)^\delta = f^\delta g + g^\delta f$.  To simplify notation,
we assume that $f,g$ are polynomials in one variable.  Let
$f(x)=\sum_{i=0}^ma_ix^i$ and  $g(x) = \sum_{j=0}^nb_jx^j$.  Then
$$
\begin{array}{lll}
(fg)^\delta&=&\sum_{k=0}^{m+n}\sum_{i=0}^k(a_i^\delta b_{k-i} + 
a_ib^\delta_{k-i})x^k\\\\
&=&\sum_{k=0}^{m+n}\sum_{i=0}^ka_i^\delta b_{k-i}x^k +
\sum_{k=0}^{m+n}\sum_{i=0}^ka_ib^\delta_{k-i}x^k\\\\
&=&f^\delta g + fg^\delta
\end{array}
$$
as desired.

Then for $(a, u) \in \tau (V)$, 
$$
\begin{array}{lll}
\tau(fg)(a,u) &=& d(fg)_a\cdot u + (fg)^\delta(a)\\\\
&=& f(a)dg_a\cdot u + g(a)df_a\cdot u + f(a)g^\delta(a)+g(a)f^\delta(a)
\\\\
&=&f(a)\tau g(a,u) + g(a)\tau f(a,u)
\end{array}
$$  
which completes the proof.
\qed

\brm
This lemma says that the map $\tau: K[V]\lra\Omt[V]$ is a {\em $\tau$-derivation},
that is, a derivation such that for all $a,b \in K$, $\dl(a)\tau(b) = \dl(b)\tau(a)$.
In fact, this map is the fundamental example of such a derivation.  
For more information, see~\cite{Ros}.
\erm

\bco
Let $V$ be a variety.  For any polynomial $F\in K[T_1,\ldots,T_m]$,
and functions $f_1,\ldots,f_m\in K[V]$,
$$
\begin{array}{lll}
\tau(F(f_1,\ldots,f_m))&=&\sum_{i=1}^m\frac{\partial F}{\partial T_i}
(f_1,\ldots,f_m)\tau f_i + F^\delta(f_1,\ldots,f_m)
\end{array}
$$
\eco
\prf
It suffices to consider $F(T_1,\ldots,T_m)=
cT_1^{a_1}\cdots T_m^{a_m}$ a monomial.  Then
$$
\begin{array}{lll}
\tau (cf_1^{a_1}\cdots f_m^{a_m}) &=&
\sum_{i=1}^m (a_i f_i^{a_i-1}\tau f_i(\prod_{j\neq i}f_j^{a_j}))+
(\tau c)(f_1^{a_1}\cdots f_m^{a_m})
\\\\
&=&\sum_{i=1}^m\frac{\partial F}{\partial T_i}
(f_1,\ldots,f_m) \tau f_i+ F^\delta(f_1,\ldots,f_m)
\end{array}
$$
as desired.
\qed

\brm      
We now define rational 1-forms and $\tau$-forms on $V$.
Notice first that the above definition actually gives a sheaf
of modules of 1-forms and $\tau$-forms.  Given any 
open set $U \sub V$, one can define $\Om [U]$ and $\Omt [U]$
as above, and for open subsets $U \sub W \sub V$,
there are natural restriction maps $\rho_{W,U}:\Om [W] \lra \Om [U]$
and $\rho_{W,U}^\tau:\Omt [W] \lra \Omt [U]$.
Define an equivalence relation on 1-forms,
where $\om_1 \in \Om [U_1]$ and $\om_2 \in \Om [U_2]$, $U_i$ non-empty,
are equivalent if they agree on $U_1 \cap U_2$ (or on any open set).
A rational 1-form is an equivalence class under this relation,
and $\Om (V)$ denotes the set of rational 1-forms.
One can then easily define the domain of a rational 1-form.
Recall that $\Om (V)$ is a $\dim(V)$-dimensional vector space
over $K(V)$, the field of rational function on $V$.
  
The rational $\tau$-forms $\Omt (V)$ are defined in exactly the same way
and also form a $K(V)$-vector space.  
As with regular \taufs, there is a natural embedding
$\iota: K(V)\lra\Omt(V)$.  
\erm

The next result follows immediately from Proposition~\ref{freemod}.

\bpro    
\label{taurat}
Given an $n$-dimensional variety $V$, $\Omt(V)$ is an 
$(n+1)$-dimensional $K(V)$-vector space.
\epro

From now on, by 1-form or $\tau$-form we will mean rational
1-form or $\tau$-form.  Following Hrushovski-Itai~\cite{HI},
regular forms will be referred to as {\em global forms}.

\brm
In the language of schemes, we can express Proposition~\ref{freemod} by
saying that, given a smooth variety $V$, $\Omt_V$ is a locally free
sheaf of dimension $\dim(V) + 1$.  In fact, elements of $\Omt_V$ are
exactly the rational functions on the prolongation cone $\tilde{V}$ that
are linear on each fiber.  Equivalently, they are sections of the dual
bundle of $\tilde{V}$.  This situation is analogous to the fact that 
1-forms are functions on the tangent bundle and sections of the cotangent
bundle.
\erm

\subs{The $\Lm$ map}
We describe a functor $\Lambda$ from $\tau$-forms to 1-forms, that
is a precise analog of the functor $\lambda$ from affine spaces to 
vector spaces, defined above.

\bla
\label{lambda}
Let $V$ be an irreducible smooth variety.  There is a natural map
$\Lambda_V: \Omega^\tau (V) \lra \Omega (V)$, such that for each 
$\omega^\tau \in \Omega^\tau(V)$, and $v\in V$,

$$
\Lambda_V(\omega^\tau)_v = \lambda(\omega^\tau_v)
$$
where $\lambda(\omega^\tau_v): TV_v\lra K$ is the linear map
associated to the affine map $\omega^\tau_v: \tau V_v \lra K$.

Given $f \in K(V)$ and $\tau f\in\Omt(V)$,
then $\Lm_V(\tau )=df\in\Om(V)$.
\ela

\prf
Fix $\omt\in\Omt(V)$, and let $\om =\Lm_V(\omt)$. By the definition of $\Lm_V$, 
for each $v\in V$, $\om_v$ is a linear map on $TV_v$, so one must show
that these linear maps vary smoothly on $V$.  It suffices to check
this locally.  

Let $v\in V$, and choose an open neighborhood $U\sub V$ of $v$, such
that on $U$, $\omt$ is given by 
$\sum_{i=1}^ng_idf_i + \iota h$, $f_i,g_i,h \in K[U]$.
Then $\om|_U=\sum_{i=1}^ng_idf_i$, as desired.

The last assertion follows from Lemma~\ref{taudiff}.
\qed
  
\msk

Notice that $\Lambda_V: \Omega^\tau(V) \lra \Omega(V)$ is clearly surjective.
The map $\Lambda_V$ restricts to a map from global $\tau$-forms to global
1-forms, but this map is in general {\em not} surjective.
Below, the remarks preceding Proposition~\ref{cohom} provide more information
about this restricted map.

Observe also that the composition of $\iota: K(V) \lra \Omt(V)$ with
$\Lambda_V: \Omt(V)$ is trivial, i.e., $\Lambda_V \circ \iota$ is the zero
map.  As $\iota$ is injective, and $\dim(\Omt(V)) = \dim(\Om(V)) + 
\dim (K(V)) = n + 1$, as $K(V)$-vector spaces, one has the following
proposition.

\bpro
For any variety $V$, there is a short exact sequence
$$
0 \lra K(V) \lra \Omt(V) \lra \Om(V) \lra 0.
$$
\epro

A description of this sequence, 
reformulated in terms of locally free sheaves on $V$
is given in~\cite{Ros}, where it is shown that, as an extension of
$\Om(V)$ by $K(V)$, it corresponds to the Kodaira-Spencer class of $V$.

We now describe how to pullback \taufs.  Recall that, given a 
rational map $\phi:V\lra W$ between varieties, there is a pullback map
$\phi^*:\Om(W)\lra\Om(V)$ that can be defined as follows.
Let $\om\in\Om(W)$, $v\in V$ and $y=\phi(v)$.  Then for all
$a\in TV_v$, $\phi^*\om_v(a)=\om_y(T\phi(a))$.

\bdf
Let $\phi:V\lra W$ be a rational map between varieties.  For each
$\omt\in\Omt(W)$, we define the {\em pullback} of $\omt$ by $\phi$, written
$\phi^{\tau*}\omt$, to be the \tauf\ on $V$, given as follows.
For $v\in V, y=\phi(v)$, and $a\in\tau V_v$, then
$\phi^{\tau*}\omt(a)= \omt_y(\phi\prl(a))$.
\edf

\bla
\label{comm}
Let $V, W$ be varieties, $\phi: V \lra W$ a morphism.  The following diagram is 
commutative.

$$
\xymatrix{
\Omt (V) \ar[d]^{\Lm_V} &\Omt (W) \ar[l]_{\phi^{\tau *}} \ar[d]^{\Lm_W}\\
\Om (V) &\Om (W)\ar[l]_{\phi^*}
}
$$
\ela

\prf It suffices to check this on fibers $\tau V_a$ and $TV_a$, $a \in V$,
$\tau W_{\phi(a)}$ and $TW_{\phi(a)}$, $\phi(a) \in W$.
We have $\phi\prl_a: \tau(V)_a \lra \tau W_{\phi(a)}$ and $T\phi_a:
TV_a \lra TW_{\phi (a)}$; notice that $\lambda(\tau\phi_{a})=d\phi_{(a)}$.
Then $(\phi^{\tau *}\omega^\tau)_a = \omega^\tau\circ\phi\prl_a$ and 
$\lambda(\phi^{\tau *}\omega^\tau)_a)=\lambda(\omega^\tau)\circ\lambda(\tau\phi_a)$.
Going around the square in the other direction, 
$\phi^*(\lambda(\omega^\tau_{\phi (a)})) = 
\lambda(\omega^\tau)\circ\lambda(\tau\phi_a)$,
as desired.
\qed
  
\brm
By the previous lemma, one can make $\Lm$ into a functor from the 
category of \taufs\ to the category of 1-forms.  Given $\phi:V\lra W$, 
we have $\Lm(\phi^{\tau*})=\phi^*$.
(Alternatively, one can define $\Lm$ explicitly and fiberwise, as in 
Lemma~\ref{lambda}.)
\erm

\sec{Algebraic curves}
\label{algcurves}


\subs{Global $\tau$-forms}
  
In this section, by curve we mean a smooth projective curve, unless noted otherwise.
Recall that $\Omega [V]$ denotes the set of global 1-forms on a variety $V$.

\bdf 
Given a variety $V$, let $\Omega^\tau [V]$ denote the set of global $\tau$-forms 
on $V$, a $K$-subspace of $\Omega^\tau (V)$.
\edf

Note that for any variety $V$, the map $\Lambda_V:\Omega^\tau(V) \lra \Omega(V)$
restricts to a map from $\Omega^\tau[V]$ to $\Omega[V]$.  The next lemma is immediate.

\bla
\label{kern}
Let $C$ be a curve.
The kernel of $\Lambda_C:\Omt [C] \lra \Om [C]$, $\textup{Ker}(\Lm_C)$, is the set of
$\omt \in \Omt[C]$ such that for all $c \in C$, $\omt_c$ is a constant affine
map on $\tau C_c$.
Thus each $\omt \in \textup{Ker}(\Lm_C)$ determines in a natural way a regular
function on $C$.  Since the only regular functions on $C$ are the 
constant functions, there is a natural isomorphism $\textup{Ker}(\Lm_C) \cong K$.

More generally, the same holds for any smooth projective variety.
\ela

We can now describe the fibers of the map $\Lm_C : \Omt [C] \lra \Om [C]$.
For all $\omt_1, \omt_2 \in \Omt [C]$, $\Lm_C (\omt_1) = \Lm_C (\omt_2)$ if and only if
$\Lm_C (\omt_1 - \omt_2) = 0 \in \Om [C]$ if and only if there is an $e \in K$ such that
for all $a \in \tau C$, $\omt_1(a) = \omt_2(a) + e$.  In this 
case, we say that $\omt_1$ is a {\em translate} of $\omt_2$, or that they
are {\em parallel}.  Clearly, parallelism classes are 1-dimensional $K$-affine
spaces and are the fibers of the map $\Lambda_C$.

Recall that on a curve $C$ of genus $g$, $\Omega[C]$ is (by definition) a 
$g$-dimensional $K$-vector space. 
Here, we note that the computation of dimension of $\Omt[C]$ is somewhat
more complicated, depending on the {\em Kodaira-Spencer class} of $C$.

\bpro
\label{globdim}
Let $C$ be a curve of genus $g$.  Then 
$1 \leq \dim_K \Omt[C] \leq g+1$.
Further, $\dim_K \Omt[C] = g + 1$ if and only if
$C$ is isomorphic to a curve defined over the field $k$ of constants.
\epro

\prf
The first statement follows immediately from the above remarks as, given the 
map $\Lm_C : \Omt [C] \lra \Om [C]$, it follows that $\dim_K \Omt [C] 
\leq \dim_K \Om[C]  + \dim_K \Ker (\Lm_C) = g +1$.  The other inequality 
follows from the natural embedding $\iota$ of $K[C] \cong K$ into $\Omt[C]$.
The proof of the second statement,which is equivalent to Proposition~\ref{cohom}, 
uses scheme-theoretic machinery developed below.
\qed

\msk

For curves of genus 1, this yields a complete description of the global $\tau$-forms.

\bco
Let $C$ be a curve of genus 1.  If $C$ is isomorphic to a curve defined over $k$,
then $\dim_K \Omt[C] =2$.  Otherwise, $\dim_K\Omt[C]=1$.  Further, in both cases,
we can give an explicit description of the global $\tau$-forms. 
\eco

\prf
Suppose first that $C$ is defined over $k$.  Then, $\tau C \cong TC$, the tangent 
variety.  In this case, it is clear that any \tauf\ can be written (uniquely) as a sum of
a 1-form and of a \tauf\ of the form $\iota f, f \in K(C)$ (as defined in 
Remark~\ref{iotaiota}), and likewise for global forms.  
Thus $\dim_K\Omt[C] = \dim_K\Om[C] + \dim_K K[C] = 1+1 = 2$.
Otherwise, if $C$ is not defined over $k$, then $\dim_K \Omt[C] = 1$, by the above
proposition.  Here, $\Omt[C] = \{\iota e | e\in K, \iota e$ a constant trivial 
$\tau$-form$\}$.  
\qed

\msk

Recall that Hrushovski and Itai call a 1-form on a curve {\em essential} if it is
not the pullback of any other 1-form, and show that there are many essential 
global 1-forms on any curve of genus $\geq 2$.  

\bpro 
\label{hrit}
\textup{(Hrushovski-Itai)}  Let $C$ be a curve of genus $\geq 2$,
defined over $k$.  There exists a finite or countable union  
${\mathcal{E}} = \cup_lS_l$ of proper subspaces of $\Om [C]$, such that any 1-form in 
$\Om [C] \backslash {\mathcal{E}}$ is essential.  
There exists an essential global 1-form, defined over $k$.
\epro

Here, we are unable to determine whether an analgous statement holds for global $\tau$-forms.
Another natural question is whether the global \taufs\ on a curve of positive genus
are exactly the pullbacks of invariant global \taufs\ on its Jacobian, as is true
of 1-forms.

\subs{Rational $\tau$-forms}

In this section, we introduce an equivalence relation on 
$\tau$-forms that will be useful for the study of strongly minimal 
sets in the next section.  Here, there is no need to assume that
curves are projective.
  
\bdf    
Let $C$ be a curve, $\pi_\tau: \tau C \ra C$ the canonical projection.
Let $\omt$ be a \tauf\, on a curve $C$.  For each $c\in C$, $\omt_c$
is either 
(i) everywhere defined on $\pi_\tau^{-1}(c)$, and a bijective map to $K$,
(ii) everywhere defined on $\pi_\tau^{-1}(c)$, but a constant map to $K$,
(iii) everywhere undefined on the affine fiber  $\pi_\tau^{-1}(c)$.
In case (ii), say that $c$ is a {\em zero} of $\omt$ (where zero here means that
given the affine map $\omt_c:\pi_\tau^{-1}(c)\lra K$, the corresponding
linear map $\lm(\omt_c): K \lra K$ is the zero map).
In case (iii), say that $c$ is a {\em pole} of $\omt$.

Let $Z_{\omt}$ denote the set of zeros of $\omt$, and
$P_{\omt}$ its set of poles.  
\edf

Observe that for any $\omt \in \Omt (C)$, either $Z_{\omt}$ is finite
or $\omt$ is trivial, in which case $Z_{\omt}$ is an open subset of $C$.

\bdf
Let $C$ be a curve, $\omt \in \Omt (C)$ non-trivial.  
The {\em null set} of $\omt$ is
$$
N_{\omt} = \{a \in \tau C | \omt(a) = 0\}.
$$
This is a rational section of the algebraic variety $\tau C$,
and thus birational to $C$.

Say that $\omt_1, \omt_2 \in \Omt (C)$ are {\em $\seq$-equivalent},
written $\omt_1\seq\omt_2$, if
$N_{\omt_1}$ and $N_{\omt_2}$ are almost equal.

Given a \tauf\ $\omt$, let $\omt / \sim$ denote its $\seq$-class.  
\edf

\bla  
\label{normform}
Let $C$ be a curve, and $\omt \in\Omt (C)$ a non-trivial $\tau$-form.
For every $\omt_1 \in \Omt(C)$, there are unique $f,g\in K(C)$ such
that $\omt_1=f\omt+\iota g$.
\ela

\prf
(The idea is similar to the proof of Proposition~\ref{freemod}.)
If $\omt_1$ is trivial, i.e., there is a $g \in K(C)$ such that
$\omt_1 = \iota g$, then $\omt_1 = 0\omt+\iota g$.  So we may
suppose that $\omt_1$ is non-trivial.

Assume first that $\omt$ and $\omt_1$ are $\seq$-equivalent.
Choose some small enough quasi-affine neighborhood of $U \sub C$
such that
\ben
\item  $U$ is disjoint from all the zeros and poles of $\omt$ and $\omt_1$;
\item  $N_{\omt}\cap U = N_{\omt_1 }\cap U$;
\item  on $U$, $\omt$ and $\omt_1$ are of the form
$$
\omt = \sum_{i=1}^mg_i\tau f_i
$$
$$
\omt_1 = \sum_{j=1}^n\gamma_j\tau \phi_j
$$
with $g_i,f_i,\gamma_j,\phi_j, \in K[U]$.
\een
In particular, $\omt$ and $\omt_1$ are regular functions on the quasi-affine 
variety $\tau U$.  Furthermore, on each fiber $\tau U_u$,
there is a single point $x \in N_{\omt} = N_{\omt_1}$ such that 
$\omt(x)=\omt_1(x)=0$.  Since $\omt_u$ and $\omt_{1,u}$ are non-constant
affine maps on $\tau U_u$ with the same zero, one must be a 
constant multiple of the other.  That is, there is an $e_u\in K$,
such that $e_u\omt_u=\omt_{1,u}$.  From the definition of $\omt$ and
$\omt_1$ it is clear that the function $e$ on $U$, $e(u) = e_u$, 
is a regular function on $U$. So $\omt_1=e\omt$, as desired.

We now consider the general case.  There is a finite set $A \sub N_{\omt}$, such that
$N_{\omt} - A$ is a regular section of $\tau U$, so there is an isomorphism 
$\psi:U\lra N_{\omt} - A$.  Let $g\in K[U]$ be the regular function
$g=\omt_1\circ \psi$, and define $\omt_2 = \omt_1 - \iota g$.
Then $\omt$ and $\omt_2$ are $\seq$-equivalent, so there is an
$f \in K[U]$ such that $\omt_2=f\omt$.  Finally,
$\omt_1 = f\omt + \iota g$, as desired.
\qed
    
\brm
Note that this lemma gives another proof that 
$\Omt(C)$ is a 2-dimensional $K(C)$-vector space,
which also follows from Proposition~\ref{taurat}.
\erm
 
\bco
\label{seqcl}
Let $C$ be a curve, $\omt_1, \omt_2 \in \Omt(C)$, both non-trivial.  
Then $\omt_1$ and $\omt_2$ are $\seq$-equivalent if and only if there is a rational 
function $f\in K(C)$ such that $\omt_1=f\omt_2$.

Given a \tauf\ $\omt$, the class $\omt / \sim$ is a 1-dimensional 
subspace of $\Omt(C)$.
\eco  

\prf
One direction of the first claim follows immediately from the above proof.  
In the other, suppose that $\omt_1=f\omt_2$.  Then on any open set $U\sub C$ 
disjoint from the zeros and poles of $f, \omt_1$, and $\omt_2$, $N_{\omt_1}|_U
=N_{\omt_2}|_U$.

The second claim follows from the the first and the observation that 
$\omt / \sim$ is closed under addition.  
\qed
      
\bla
\label{sigmanull}
Let $C$ be a curve, $\omt_1, \omt_2 \in \Omt(C)$, both non-trivial.  
Suppose that $N_{\omt_1} \cap N_{\omt_2}$ is infinite.  Then 
$\omt_1$ and $\omt_2$ are $\seq$-equivalent.  In other words,
if  $N_{\omt_1} \cap N_{\omt_2}$ is infinite, then 
$N_{\omt_1}$ and $N_{\omt_2}$ are almost equal.
\ela
  
\prf
$N_{\omt_1}$ and $N_{\omt_2}$ are each rational
sections of $\tau C$, so each is a curve on the surface $\tau C$.
If these curves have infinite intersection, they must be equal, up 
to a finite set.
\qed

\msk

Since we will be concerned with $\sim$-equivalence classes of 
\taufs, the following relativized versions of global and essential \taufs\ 
are important for the main results about strictly minimal sets below.  
Equivalent, perhaps more natural, definitions of 
these notions, in the language of schemes, are given in Section~\ref{tauinv}.

\bdf  Let $C$ be a smooth projective curve, $\omt$ a rational \tauf.
Say that $\omt$ is $\sim$-{\em global} if it is $\sim$-equivalent to 
a global \tauf, in which case we also say that the class $\omt / \sim$ is 
{\em global}.
Say that $\omt$ is $\sim$-{\em essential} if every
$\sim$-equivalent \tauf\ is essential.  
\edf  
 
\brm
Let $C, C'$ be curves, and $f:C \lra C'$ a morphism between them.
It is easy to see that for any $\sim$-equivalent \taufs\ $\omt_1$
and $\omt_2$ on $C'$, then  $f^*\omt_1$ and $f^*\omt_2$ 
are $\sim$-equivalent on $C$.
\erm

\brm
Let $\omt$ be a non-trivial global \tauf\ on a curve $C$ of genus $g \geq 1$.  
Then for any $c \in K$, $\omt_c : = c + \omt$ is global and not 
$\sim$-equivalent to $\omt$.  In  particular, there is a family
of global $\sim$-classes of $\tau$-forms indexed by $K$.  
\erm

Recall that we defined two global \taufs\ on a curve $C$ to be parallel 
if they are both in the same fiber of the map $\Lm_C$.  This notion also
makes sense for rational \taufs.

\bdf
Let $C$ be a curve.  We say that two \taufs\ $\omt_1$ and $\omt_2$ are 
{\em parallel}, written $\omt_1\teqq\omt_2$, if there is a $g\in K(C)$ such that 
$\omt_2=\omt_1+g$.  Equivalently, $\Lm_C(\omt_1) = \Lm_C(\omt_2)$.
\edf

Parallelism is clearly an equivalence relation.

\bla  Fix a curve $C$.
\ben
\item  The trivial \taufs\ on $C$ are a parallelism class.
\item  On the non-trivial \taufs, $\seq$-equivalence and parallelism 
  are cross-cutting equivalence relations.  More precisely, for 
  any non-trivial $\omt_1, \omt_2$,
  $|\omt_1/ \seq \cap \omt_2/ \teqq|=1$.
\een
\ela
\prf
Immediate.
\qed

\para{The $\Lm$ map on curves}  
Fix a curve $C$, and a non-trivial \tauf\ $\omt_0\in\Omt(C)$.  
By Lemma ~\ref{normform}, every $\omt\in\Omt(C)$ 
is equal to $f\omt_0+g, f,g\in K(C)$, and, by Corollary~\ref{seqcl}
the $\seq$-equivalence class of $\omt_0$ is 
$\omt_0/ \seq=\{f\omt_0:f\in K(C)-\{0\}\}$.

The following lemma can be easily verified.

\bla
Let $C$ be a curve, $f\in K(C)$, and $\omt$ a non-trivial $\tau$-form on $C$.  
Then $\Lm_C(f\omt)=f\Lm_C(\omt)$.
%
Thus the map $\Lm_C$ induces a bijection between the 
equivalence class $\omt/ \seq$ and the non-trivial 1-forms.
\ignore{
\item  For any $\omt_1,\omt_2\in\Omt(C)$, $\Lm_C(\omt_1)=\Lm_C(\omt_2)$ if 
and only if $\omt_1\teqq\omt_2$.  Equivalently, for any $\om\in\Om(C)$, the 
fiber $\Lm^{-1}_C(\om)$ is a parallelism class.  
}
\ela

\sec{Order one strongly minimal sets}
\label{sm}  

We aim to classify order one strongly minimal sets, up to
non-orthogonality.  
Hrushovski proved that every such set is either
trivial or non-orthogonal to the constant field $k$, so we will
be investigating trivial sets.
We first briefly recall the analysis from Hrushovski-Itai of order 
one \sm sets over the constants.
  
Say that two strongly minimal sets $X$ and $Y$ are {\em birationally 
isomorphic}, or {\em birational}, if there is an almost everywhere defined 
bijective map from $X$ to $Y$.  Any order one strongly minimal set $Y$ is 
birational to one of the form
$$
X = \Xi(C, s) = \{ x \in C : \dl x = s(x)\}
$$
where $s: C \lra \tau (C)$ is a rational section of the prolongation.
In this case, one says that $X$ lives on $C$, and the differential
function field of the Kolchin closed set $X$ equals the function 
field of $C$.

Suppose now that $C$ is a curve defined over $k$, so $\tau C = TC$,
the tangent space of $C$.  Then the \sm sets $X$ defined over $k$
and living on $C$ are given as $X = \Xi (C, s)$, $s$ a 
rational vector field also defined over $k$.  There is a natural
1-1 correspondence between such vector fields $s$  and 1-forms
$\om$, given by the relation $\om(x)s(x) = 1$ for almost all $x$.  
So order one \sm sets over the constants also correspond to sets
$$
X = \Xi(C, \om) = \{a \in C : \om(a)\dl a=1\}
$$
where this is taken to include a pole $p$ of $\om$ if $\dl p = 0$.
(This last condition ensures that $\Xi(C, \om)$ is Kolchin closed.)

In this setting, Hrushovski and Itai analyze \sm sets in terms of 
pairs $(C, \om)$, $C$ a smooth projective curve and $\om \in \Om (C)$,
both defined over $k$.  We now show how to extend their results
to all curves, using $\tau$-forms.

\bdf
Let $C$ be a smooth projective curve and  $\omt \in \Omt (C)$ a
non-zero $\tau$-form.  Let $P \sub C$ be the set of poles $p$ of $\omt$
such that $\dl p=0$, and define
$$
\Xi (C, \omt) = \{a \in C : \omt(a)\dl a = 0\} \cup P
$$
\edf
As above, one can check that $\Xi (C, \omt)$ is always Kolchin closed.

Observe that $X = (C, \omt) = \{a \in C : \omt(a)\dl a = 0\} $ is a strongly minimal 
set.  Indeed, let $s$ be a section of $\tau C$ so that $\omt s = 1$ almost 
everywhere.  Then $\Xi(C, s) = \Xi(C, \omt)$, up to a finite set.
On the other hand, 
this is no longer a 1-1 correspondence, as many $\tau$-forms correspond to the 
same section and therefore determine the same \sm set.  

\bla
Let $C$ be a curve, $\omt_1, \omt_2 \in \Omt(C)$, both non-trivial.
Then $\omt_1$ and $\omt_2$ are $\seq$-equivalent if and only if 
$\Xi(C, \omt_1)$ and $\Xi(C, \omt_2)$ are almost equal.
\ela

\prf
From the definition, $\Xi(C,\omt)$ only depends on $N_{\omt}$,
so if $\omt_1$ and $\omt_2$ are $\seq$-equivalent, then 
$\Xi(C, \omt_1)$ and $\Xi(C, \omt_2)$ are almost equal.

In the other direction, if
$\Xi(C, \omt_1)$ and $\Xi(C, \omt_2)$ are almost equal,
then $N_{\omt_1}\cap N_{\omt_2}$ is infinite, so 
$\omt_1$ and $\omt_2$ are $\seq$-equivalent, by 
Lemma~\ref{seqcl}.
\qed

\bco
\label{formz}
There is a 1-1 correspondence between $\seq$-equivalence
classes of \taufs\, in $\Omt (C)$ and strongly minimal sets living on $C$.
\eco
Notice that we do not have to treat the case $C(k) \sub C$, corresponding
to $s(x)=0$ separately.  
   
The following series of lemmas provides versions of 
Lemmas 2.4 and 2.6 -- 2.9 of \cite{HI}, for \taufs.  
(Lemma 2.5 becomes false.)  The proofs are the same.

For the proof of the next lemma, we recall the following easy fact.
Given $\phi:V\lra W$, and $v\in V$, $(\tau \phi)(\dl v)=\dl (\phi(v))$.

\bla  
\label{taupullback}
Let $g: C_1 \lra C_2$ be a dominant  regular map between nonsingular curves, and 
$\omt_2$ a \tauf\, on $C_2$.  Let $\omt_1 = g^{\tau *}\omt_2$ be the pullback of
$\omt_2$ by $g$ to $C_1$.  Then 
$g^{-1}\Xi(C_2, \omt_2) = \Xi(C_1, \omt_1)$.  
\ela

\prf  
Let $c_1\in C_1, c_2=g(c_1)$.  Then
$$
\omt_1(c_1)\dl c_1=\omt_2(c_2)(g\prl)(\dl c_1)=\omt_2(c_2)\dl 
(g(c_1))=\omt_2(c_2)\dl c_2.
$$
Thus $\omt_1(c_1)\dl c_1=0$ if and only if $\omt_2(c_2)\dl c_2=0$.
\qed

\bla
\label{infint}
If $\Xi(C, \omt_1)$ and $\Xi(C, \omt_2)$ have infinite intersection,
$\omt_1 $ and $\omt_2$ are $\seq$-equivalent.
\ela

\prf
For infinitely many points $c$ in the intersection $\Xi(C, \omt_1)
\cap \Xi(C, \omt_2)$, $\omt_1(c)\dl c=0=\omt_2(c)\dl c$.  Thus,
$N_{\omt_1}$ and $N_{\omt_2}$ have infinite intersection, and are
therefore $\seq$-equivalent, by Lemma~\ref{sigmanull}.
\qed

\bla
\label{compat}
Let $g:C_1 \lra C_2$ be a dominant rational map between non-singular curves.
Let $\omt_i$ be a \tauf\, on $C_i, i = 1,2$.  Assume
$g\Xi(C_1, \omt_1) \sub \Xi(C_2, \omt_2)$ (perhaps up to a finite set).
Then $\omt_1$ and $g^{\tau *}\omt_2$ are $\seq$-equivalent.
\ela  
   
\prf
By Lemma~\ref{taupullback}, $g^{-1}\Xi(C_2, \omt_2) = \Xi(C_1, g^{\tau *}\omt_2)$
(a.e.).  Thus the intersection of $\Xi(C_1, \omt_1)$ with $\Xi(C_1, g^{\tau *}\omt_2)$
is infinite.  By the previous lemma, $\omt_1$ and $g^{\tau*}\omt_2$
are $\seq$-equivalent.  
\qed

\bla
Let $C_1$ be a complete nonsingular curve, $\omt_1$ a global \tauf\, on $C_1$.
Let $g:C_1 \lra C_2$ be a rational function into a complete nonsingular curve,
and $\omt_2$ a \tauf\, on $C_2$ such that $g\Xi(C_1, \omt_1) = 
\Xi(C_2, \omt_2)$.  Then $\omt_1$ is a $\sim$-global \tauf\ on $C_2$.
\ela
  
\prf
As pointed out in \cite{HI}, $g: C_1 \lra C_2$ is a surjective morphism.
It then follows from definitions that $\omt_2$ is a global \tauf\, on $C_2$
if and only if its pullback, $g^{\tau *}\omt_2$ is a global \tauf\, on $C_1$.
By Lemma~\ref{taupullback}, $g^{-1}\Xi(C_2, \omt_2)=\Xi(C_1, g^{\tau *}\omt_1)$,
so $\Xi(C_1, g^{\tau *}\omt_2) \cap \Xi(C_1, \omt_1)$ is infinite.
By Lemma~\ref{infint}, $g^{\tau *}\omt_2 \sim \omt_1$, so by
the above remark, $\omt_1$ is a $\sim$-global form.
\qed  

\bla
\label{noeqrel}
Let $\omt$ be a $\sim$-essential \tauf\ on a curve $C$.  Let $E$ be a 
definable equivalence relation on $\Xi(C,\omt)$, with finite classes.
Then almost every class of $E$ has one element.
\ela

\prf
Suppose not.  By \cite{HI}, Lemma 1.2, there is a curve $C_0$ and a map
$g: C \ra C_0$ with deg$(g) > 1$ such that $g(\Xi(C,\omt))$ lives on $C_0$.
Thus there is a \tauf\ $\omt_0$ on $C_0$ such that $g(\Xi(C,\omt))=
\Xi(C_0,\omt_0)$, up to a finite set.  But by Lemma~\ref{compat},
$\omt$ and $g^{\tau*}(\omt_0)$ are $\sim$-equivalent, contradicting the
assumption that $\omt$ is $\sim$-essential.
\qed

\bla
\label{nonorthcon}
Let $C$ be a curve of genus $\geq 1$ and $\omt$ a \tauf\ on $C$.  
If $\Xi(C,\omt)$ is non-orthogonal to the constants, then there is a 
regular map $g:C \ra \P^1$ with $g(\Xi(C,\omt))=k$, up to a finite set.  
In particular, there is a \tauf\ $\omt_0$ on $\P^1$ such that 
$\Xi(\P^1,\omt_0) = k$ and $g^{\tau*}\omt_0$ is $\sim$-equivalent 
to $\omt$ on $C$.  Thus, $\omt$ is not $\sim$-essential.
\ela

\prf
Let $X = \Xi(C,\omt)$ be non-orthogonal to the constants.  Then there is 
a definable differential rational function $f$ with $f(X) \sub k$, up to a
finite set (e.g., see \cite{HI}, proof of Lemma 2.10).  Since for all 
$x \in X, \dl(x)=s(x)$, for some polynomial $s$, we can assume that
$f$ is a rational function.  Thus $f$ extends to a regular, dominant map 
$f:C \ra \P^1$, proving the first claim.  The rest of the lemma follows
from Lemma~\ref{compat}. 
\qed

\msk

We are now able to generalize a central result of \cite{HI}.  Recall
that a strongly minimal set $X$ defined over a set $A$ is {\em strictly minimal}
over $A$ if $X \cap \acl(A) = \emptyset$ and, for all $a \in X, \acl(aA) \cap X
= \{a\}$.  We say that $X$ is strictly minimal if it is strictly minimal over
some set over which it is defined.  A strongly minimal set $X$ 
is {\em $\omega$-categorical} if for any finite set $B \sub X$, $\acl_X(B)$
is finite.  Hrushovski's theorem that all order one strongly minimal\ sets are 
$\omega$-categorical~\cite{HJ} plays an important role in the proof of the 
following result of \cite{HI}.

\bpro
\label{strictmin}
Let $C$ be a complete nonsingular curve over $k$ of genus $>1$.  Let $\om$ 
be an essential 1-form on $C$, defined over $k$.  Then, after perhaps removing
a finite set, $\Xi(C,\om)$ is strictly minimal, with trivial induced 
structure.  Two such sets $\Xi(C_1,\om_1)$
and $\Xi(C_2,\om_2)$ are non-orthogonal if and only if there exists an isomorphism 
$g:C_1 \ra C_2$ with $\om_1 = g^*\om_2$.
\epro
(In \cite{HI}, they state this result for global 1-forms, which implies
additionally that $\Xi(C,\om)$ has no points in the algebraic closure of
the parameters.)

We now state one of our main results.

\bthm
\label{genprop}
Let $C$ be a smooth projective curve of genus $g \geq 1$.  Let $\omt$ be
a $\sim$-essential \tauf.  Then there is a finite set $A \sub \Xi(C,\omt)$
such that $\Xi(C,\omt) - A$ is trivial strictly minimal. Further, for any
definable equivalence relation with finite classes on $\Xi(C,\omt)$, 
over any set of parameters, almost every class has one element.

Two such sets, $\Xi(C_1,\omt_1)$ and $\Xi(C_2,\omt_2)$ are 
not orthogonal if and only if there is an isomorphism $f:C_1 \ra C_2$ such that
$f^{\tau *}\omt_2$ is $\sim$-equivalent to $\omt_1$.
\ethm

\prf
By Lemma~\ref{nonorthcon}, $X = \Xi(C,\omt)$ is a trivial strongly minimal set and
by Lemma~\ref{noeqrel}, for any definable equivalence relation with finite classes 
on $X$, over any set of parameters, almost every class has one element.
Since $X$ is $\omega$-categorical, the set $A$ of elements of $X$ which are 
algebraic over the parameters defining $C$ and $\omt$ is finite.  Thus 
$X - A$ is strictly minimal.

Suppose $\Xi(C_1,\omt_1)$ and $\Xi(C_2,\omt_2)$ are non-orthogonal.  
After removing a finite set of points from each, one is left with two
strictly minimal sets with trivial geometry, so there is a definable 
bijection between them.  Thus there is a birational map $g: C_1 \ra C_2$,
with $\omt_1 \sim g^{\tau *}\omt_2$, which is an isomorphism.
\qed

\brm
In a sense, this completely classifies the trivial strictly minimal sets,
in terms of $\sim$-essential forms, but it is not very explicit.  
In particular, it leaves open the question whether there are any such forms.
For curves over $k$, one gets an affirmative answer from \cite{HI}, but the general
case remains open.

A key part of Hrushovski and Itai's proof that there are many trivial strictly 
minimal sets living on every curve defined over $k$ of genus $>1$,
which are all orthogonal, is the fact that on such curves generic global 
1-forms are essential.  Proving an analogous statement for \taufs\ seems difficult,
since the dimension of global \taufs\ does not depend only on the genus.
Further, in the more general case, one needs to consider $\sim$-essential \taufs,
which are essentially equivalence classes.  For these reasons, it is not
clear how one might be able to adapt ideas from \cite{HI} to answer this question.
\erm

Let us call an equivalence relation on a set $X$ {\em trivial} if it is definable
in the language of equality.  Equivalently, there is either a cofinite equivalence class
or all but finitely many elements are in classes of size 1.

In the statement of Proposition~\ref{strictmin}, there are no non-trivial equivalence 
relations on $\Xi(C,\omega)$ defined over $k$.   One might ask whether it is possible
to define such a relation using parameters from $K$.  If this were the case, then
one would have pairs $(C,\omega)$ and $(C_0,\omt_0)$, $C \not\cong C_0$,
$C$ a curve defined over $k$, $\omega$ an essential
1-form defined over $k$, $\Xi(C,\omega)$ trivial strictly minimal over $k$, and
$\Xi(C_0,\omt_0)$ strictly minimal over $K$, with  a finite-to-one map from 
$\Xi(C,\omega)$ to $\Xi(C_0,\omt_0)$.    Below, we argue that this can not in fact
happen.  The key point is the following lemma.
  
\bla
\label{eqprm}
Let $K$ be an $\om$-stable structure, and let $X$ be a trivial strictly minimal set in $K$
defined over a set $A$.  Then no non-trivial equivalence relation on $X$ is definable
with parameters from $K$.
\ela

\prf
As usual, we may assume that the set $A$ over which $X$ is defined is empty.
Suppose for contradiction that there is a tuple $\abar$ and a formula $\phi(\xbar,y,z)$
such that $\phi(\abar,y,z)$ defines a non-trivial equivalence relation $E$ on $X$.
By strong minimality of $X$, it is clear that there is a finite $n$ such that every
class in $E$ is of size $\leq n$ and all but finitely many classes are of the same
size $m > 1$.  Let $\theta(\xbar)$ be a formula true of $\abar$ such that for all
$\bbar$, if $K \models \theta(\bbar)$, then $\phi(\bbar,y,z)$ defines an equivalence relation
on $X$ such that all but finitely many elements of $X$ belong to a class of size $m$.

Let $T$ be the set of equivalence relations $T = \{\phi(\bbar,y,z) | K \models \theta(\bbar)\}$.
Note that, even without elimination of imaginaries in $K$, one can treat $T$ as a definable
set that can be quantified over in a natural way.  Observe also that $T$ is infinite
as follows.  For if not, choose an $E \in T$ and a pair $(p,q)$ in $E$.  
Then $q \in \acl (p)$, contradicting the fact that $X$ is trivial strictly minimal.

Let $\Sigma$ be the set of finite sequences of 0's and 1's, partially ordered such that,
for $\rho, \sigma \in \Sigma$,
$\rho \leq \sigma$ if and only if $\rho$ is an initial segment of $\sigma$.
For each $\sigma \in \Sigma$, we define a pair $\pi_\sigma = (p_\sigma,q_\sigma)$
of elements in $X$ and a `definable' set $T_\sigma \sub T$ with the following properties.  
\ben
\item $T_\sigma$ is infinite and definable over parameters 
$\cup_{\rho < \sigma}\{p_\rho,q_\rho\}$.

\item $T_\sigma$ is the disjoint union of $T_{\sigma,0}$ and $T_{\sigma,1}$.

\item $T_{\sigma,0}:=\{E \in T_\sigma | E(p_\sigma,q_\sigma)\}$ and 
      $T_{\sigma,1}:=\{E \in T_\sigma | \neg E(p_\sigma,q_\sigma)\}$ are both non-empty.
\een
Since it is clear that such a set of $T_\sigma$'s contradicts $\omega$-stability,
this will yield the desired contradiction.

We proceed by induction.  
Let $T_\emptyset = T$.  Choose a pair $\pi_\emptyset=(p_\emptyset,q_\emptyset)$ such that 
there are $E,E' \in T_\emptyset$ 
with $E(p_\emptyset,q_\emptyset)$ and $\neg E'(p_\emptyset,q_\emptyset)$.
(If this were not possible, then $T$ would only contain one element.)  

Suppose now that for all $\rho \leq \sigma$, $\pi_\rho$ and $T_\rho$ have been defined.
Let $T_{\sigma,0} = \{E \in T_\sigma| E(p_\sigma,q_\sigma)\}$ and 
$T_{\sigma,1} = \{E \in T_\sigma| \neg E(p_\sigma,q_\sigma)\}$.
By construction, both $T_{\sigma,0}$ and $T_{\sigma,1}$ are non-empty.
In fact, we claim that both sets must be infinite.  For if not, it is easy to see
that $q_\sigma \in \acl(\{p_\rho,q_\rho | \rho < \sigma\}\cup\{p_\sigma\})$,
contradicting the fact that $X$ is trivial strictly minimal.
Now choose $\pi_{\sigma,i}=(p_{\sigma,i},q_{\sigma,i})$, $i = 0,1$,  disjoint from all
$\pi_\rho$, $\rho \leq \sigma$, such that 
there are $E_i, F_i \in T_{\sigma,i}$ such that $E_i(p_{\sigma,i},q_{\sigma,i})$
and $\neg F_i(p_{\sigma,i},q_{\sigma,i})$.

It is clear that this construction has the required properties.
\qed

\bco
Let $C$ be a complete nonsingular curve over $k$ of genus $>1$.  Let $\om$ 
be an essential 1-form on $C$, defined over $k$, so that after perhaps removing 
a finite set, $\Xi(C,\om)$ is strictly minimal with trivial induced 
structure.  Then there is no non-trivial equivalence relation definable on $\Xi(C,\om)$,
even with arbitrary parameters.
\eco
  
\prf
By Proposition~\ref{strictmin} and Lemma~\ref{eqprm}.
\qed

\section{Another description of $\tau$-forms}
\label{tauinv}


In this section, we give a new algebraic characterization of the 
strongly minimal sets living on a curve $C$, using the theory of $\tau$-differentials
developed in \cite{Ros}.   Of course, it is equivalent
to the earlier formulation, but has the advantage that it does not involve
considering equivalence classes, so it is easier to work with.
Below, we use the language of schemes, so a variety is an integral separated
scheme of finite type over an algebraically closed field.
We continue to assume that curves are smooth and projective.

Recall that there is a natural bijection between rank $n$ locally free
sheaves on a variety and rank $n$ vector bundles over it.
For example, the sheaf of differential forms on a variety $X$,
$\Om_X$, corresponds to the cotangent bundle $T^*X$, 
and can be viewed as the sections of this bundle.  Given a curve $C$,  
the locally free rank 2 sheaf of $\tau$-differentials on $C$, as defined in \cite{Ros},
and also denoted $\Omt_C$, 
corresponds exactly to the set of rational $\tau$-forms, as defined above,
with a $\tau$-form $\omt$ corresponding to a section $s \in \Omt(U)$,
for open $U \sub C$.

Given a rank 2 bundle, one can consider all the rank 1 subbundles.
(To be precise, a rank 1 subbundle is a rank 1 bundle together
with a bundle embedding.)  Passing to sheaves, rank
1 subbundles correspond to invertible (i.e., rank 1) subsheaves.
Note that it is possible for one rank 1 sheaf to be properly embedded in
another (in which case the quotient will be a torsion sheaf), essentially
because this can also happen with free rank 1 modules.  Nevertheless,
in our situation, any invertible subsheaf of a rank 2 locally free sheaf 
embeds in a unique maximal invertible subsheaf, as described below.

We have seen that \taufs\ correspond to sections of a rank 2 locally
free sheaf.  Since the $\seq$-equivalence class of such a form $\omt$
consists of those forms $g\omt, g \in K(C)$, such a class should correspond
to a (unique) maximal subsheaf containing $\omt$ as a section on an open set.
The following easy lemma makes this precise.

\bla
Let $C$ be a smooth curve, $\MH$ a rank 2 locally free sheaf on $C$,
and $s \in \MH(U)$ a section.  Then there is a unique maximal 
invertible subsheaf $\MF \sub \MH$ such that $s \in \MF(U)$.  
Further, the quotient $\MH / \MF$ is an invertible sheaf.
\ela

\prf
By possibly passing to a smaller open set, we may assume that $U = \Spec\ R$
is affine, $\MH(U) \cong \tilde{M}$, $M = \<m,n \>$ a free rank 2 
$R$-module, and $s = am + bn\in M, a,b \in R$.  Since $R$ is a UFD 
(e.g.\ \cite{Mat}, Thm.\ 20.1), $a$ and $b$ have a greatest common divisor
$c$.  Letting $a' = a/c$, $b'=b/c$, and $s'=a'm + b'n$, the module 
$N = Rs' \sub M$ is the unique largest rank 1 free submodule of $M$
containing $s$.  Let $d,e \in R$ be such that $da'+eb'=1$, and let
$L = R(em-dn)$.  Then it is easy to check that $M  = N \op L$, and hence that 
there is a short exact sequence $0 \lra N \lra M \lra L \lra 0$.  

We want to define $\MF$ such that $\MF(U) = \tilde{N}$.  Since a sheaf is 
determined by its stalk at each point $P$, it suffices to show how, for each 
$P \in C - U$, $\MF(U)$ determines a unique maximal stalk $\MF_P$
that will make $\MF$ locally free.  

For each such $P$, choose an affine scheme $\ W \sub U \cup \{P\}$
containing $P$, $W \cong \Spec\ S$, such that $\MH(W)$ is free on $W$.  
Note that $ W':= W - P =
\Spec\ S_f$, for some $f \in S$, and $P = V(f)$.  Let $\MH(W') = \tilde{J}$
and choose basis elements $\<i,j\> \in J$ such that $\MF(W') = \tilde{I}$,
where $I := Ri \sub J$.  We then define $\MF_P$ to be the submodule of
$\tilde{J}_P$ generated by $\<i\>$.  
Clearly, this determines an extension of $\MF$ to $U \cup P$ and thus,
continuing in this way to all of $C$.  (As an aside, note that one could 
also define $\MF_P$ to be generated by $f^mi$, for any positive $m$.
Choosing $m > 1$ would yield a non-maximal invertible subsheaf.)
\qed

\msk

It is easy to see that two sections determine the same invertible subsheaf
if and only if, on some non-empty open set on which they are both defined,
each is a multiple of the other.  In particular, $\sim$ equivalence classes
of \taufs\ correspond to maximal invertible subsheaves, as desired.

\bdf
Let $C$ be a curve.  A {\em $\tau$-invertible sheaf} on $C$, denoted $\MFT$,
is a maximal invertible subsheaf of $\Omt_C$, i.e, an invertible sheaf $\MF$, 
together with an embedding $f:\MF\ra\Omt_C$.
\edf

\brm
For any curve $C$, there is a {\em trivial} \tis, corresponding to
the natural embedding of $\O_C$ into $\Omt_C$.  Clearly, this \tis\
corresponds to the $\sim$-equivalence class of trivial \taufs.  
Below, by \tis\ we will always mean non-trivial \tis.
\erm

We can reformulate Corollary~\ref{formz} as follows.
  
\bpro
Let $C$ be a curve.  There is a natural bijection between strongly 
minimal sets living on $C$ and $\tau$-invertible sheaves on $C$.
\epro

As before, given a \tis\ $\MFT$ on $C$, let $\Xi(C,\MFT)$ denote the
corresponding strongly minimal\ set.


We now introduce the analog of a global \tauf.  
Since we are dealing with equivalence classes of \taufs,
the correct notion is that of being a \tauf\ that is $\seq$-equivalent
to a global \tauf.  Passing to $\tau$-invertible sheaves, these are
the sheaves $\MFT$ that contain a global section.

\bdf
Let $C$ be a curve, $\MFT$ a \tis. Say that $\MFT$ is a {\em global \tis}
if $H^0(C,\MFT) \neq 0$, that is, if $\MFT$ has a global section.
\edf

Next, we define the appropriate notion of a $\tau$-invertible
sheaf being essential, which is more subtle than for 1-forms.
First, we establish the following results.

\bpro
\label{injpull}
Let $f: X \ra Y$ be a finite morphism of curves.  
There is an exact sequence of sheaves on $X$,
$$
0 \lra f^*\Omt_Y \sr{\al}{\lra} \Omt_X \lra \Om_{X/Y} \lra 0
$$
\epro

\prf
To check that $\alpha$ is injective, it suffice to check this at the 
stalk of the generic point.  Here, injectivity is obvious.
\qed

The following proposition follows immediately.

\bpro
\label{tauinvpull}
Let $f:X \ra Y$ be a finite morphism of curves, and $\MFT$ a (non-trivial)
\tis\ on $Y$.
Then there is a unique  \tis\ $\MGT$ on $Y$ such that $\al(f^*\MFT) \sub \MGT$.
Further, there is the following commutative diagram, with exact rows,
$$
\xymatrix{
0 \ar[r] & f^*\MFT \ar[r]\ar[d] & \MGT \ar[r]\ar[d] & \Om_{X/Y} \ar[r]\ar[d] & 0 \\
0 \ar[r] & f^*\Omt_Y \ar[r]^\al & \Omt_X \ar[r] & \Om_{X/Y} \ar[r] & 0,\\
}
$$
such that each vertical arrow is an injection, and $\MH$ is a torsion sheaf.
\epro

\bpr
The only point that needs to be checked is that the cokernel of the map
$f^*\MFT \lra \MGT$ is $\Om_{X/Y}$.  To verify this, it suffices to compute
the map locally around each ramification point of $X$.
\epr

\brm
This proposition says that, given a morphism of curves,
$f:X\ra Y$, the pullback of a \tis\ is not necessarily
a \tis, even though the pullback of a $\tau$-form is a $\tau$-form.
By analogy, the pullback of a 1-form on $Y$ is a 1-form on $X$
though, in general,  $f^*\Om_Y \neq \Om_X$.  But there is an 
injective homomorphism $f^*\Om_Y \ra \Om_X$, just as there is
an injective homomorphism $f^*\MFT \ra \MGT$.

Note that given a morphism between curves (or more generally varieties), 
$f: C_1 \lra C_2$, the pullback of the trivial \tis\ on $C_2$ is equal to 
the trivial \tis\ on $C_1$, precisely because $f^*\O_{C_2} = \O_{C_1}$.
\erm

\bdf
In the notation of the previous proposition, we say that $\MGT$ is the 
{\em weak pullback} of $\MFT$ (along $f$).  
In case the natural map $f^*\MFT \ra \MGT$ is an isomorphism, we say 
that $\MGT$ is the {\em strong pullback} of $\MFT$.
\edf

The following lemma follows immediately from the definitions.

\bla
Let $f: X \ra Y$ be a morphism of curves, $\omt$ a $\tau$-form on $Y$,
and $\MFT$ the \tis\ corresponding to $\omt / \sim$.  Then $\MGT$, the weak
pullback of $\MFT$, is the \tis\ corresponding to $f^*\omt$.
\ela

We can now reformulate Lemma~\ref{taupullback} in terms of \tiss.
\bla
Let $f: X \ra Y$ be a morphism of curves, $\MFT$ a \tis\ on $Y$,
and $\MGT$ the weak pullback of $\MFT$ on $X$.  Then
$f^{-1}\Xi(Y,\MFT)=\Xi(X,\MGT)$.
\ela

\bdf
Let $X$ be a curve.  A \tis\ $\MGT$ on $X$ is {\em essential}
if there does not exist a curve $Y$ and a morphism $f:X\ra Y$ 
such that $\MGT$ is the weak pullback of some \tis\ on $Y$.
\edf

\brm
Let $X$ be a curve, and $f:X \ra Y$ a regular map to another curve $Y$.
By Proposition~\ref{tauinvpull}, there is a natural embedding of $\O_X$-sheaves,
$\alpha: f^*\Omt_Y \lra \Omt_X$.  Notice that for any \tis\ 
$\MGT\sub \Omt_X$, there is maximal invertible subsheaf of 
$\alpha(f^*\Omt_Y)$ that embeds in $\MGT$, 
namely $\MGT \cap \alpha(f^*\Omt_Y)$.
Let us call this subsheaf $\MH$.  If $\MGT$ is essential, then
$\MH$ cannot be of the form $f^*\MFT$, for $\MFT$ a \tis\ on $Y$.

In particular, although every section of $f^*\Omt_Y$ is naturally a
function on $\tau X$, not every such function is the pullback along
$f$ of a $\tau$-form on $Y$.  Functions arising in this way
naturally correspond instead to sections of $f^{-1}\Omt_Y$, but this 
is not an $\O_X$-sheaf.  Recall that by definition, 
$f^*\Omt_Y = f^{-1}\Omt_Y \ot_{f^{-1}\O_Y}\O_X$.

As sections of a \tis\ $\MFT$ on $Y$ correspond to a $\sim$-equivalence class of functions 
on $\tau Y$, sections of $f^{-1}\MFT$ correspond to the pullback of these
functions to $\tau X$ along the lifting map $f\prl$, as defined in Section~\ref{prolong}.
As $f^*\MFT = f^{-1}\MFT\ot_{f^{-1}\O_Y}\O_X$, sections of $f^*\MFT$ correspond to 
products of functions from $f^{-1}\MFT$ with rational functions on $X$.
\erm

\brm
Observe that a $\tau$-form $\omt$ is $\sim$-essential, as defined in
Section~\ref{algcurves}, 
if and only if the \tis\ corresponding to $\omt / \sim$ is essential.
\erm

We can now restate Theorem~\ref{genprop}.

\bthm
Let $C$ be a smooth projective curve of genus $g \geq 1$.  Let $\MFT$ be
an essential \tis.  Then there is a finite set $A \sub \Xi(C,\MFT)$
such that $\Xi(C,\MFT) - A$ is strictly minimal. Further, for any
definable equivalence relation with finite classes on $\Xi(C,\MFT)$, 
over any set of parameters, almost every class has one element.
Two such sets, $\Xi(C_1,\mathcal{F}^\tau_1)$ and 
$\Xi(C_2,\mathcal{F}^\tau_2)$ are 
non-orthogonal if and only if there is an isomorphism $f:C_1 \ra C_2$ such that
$\mathcal{F}^\tau_1$ is the weak pullback of $f^{*}\mathcal{F}^\tau_2$.
\ethm

We now complete the proof of Proposition~\ref{globdim}.  Given a (locally free) 
sheaf $\MF$ on a projective variety $V$, let $\hh^0(V,\MF)= \dim_K H^0(V,\MF)$,  
the dimension of the space of global sections of $\MF$.

Given a curve $C$, we will call $\hh^0(C,\Omt_C)$ the {\em $\tau$-rank} of $C$, 
denoted $g^\tau_C$ or simply $g^\tau$.  The following proposition shows that $g^\tau$
can be easily determined by the Kodaira-Spencer rank of $C$, which we now explain.

By \cite{Ros}, we have the following short exact sequence 
$$
\xi:\qquad \ 0 \lra \O_C \lra \Omt_C \lra \Om_{C/K} \lra 0.
$$
Recall that $\xi$ is called an extension of $\Om_{C/K}$ by $\O_C$,
and that there is a natural bijection between the set of such extensions and 
the cohomology class $\Ext^1(\Om_{C/K}, \O_C)$, such that 
the split extension corresponds to $0\in \Ext^1(\Om_{C/K}, \O_C)$ 
(e.g., \cite{CW}, Theorem 3.4.3).  In our case, this mapping is given as follows.
Given $f:C \lra \Spec\, K$, consider the sequence 
$$
\eta: \qquad 0 \lra f^*\Om_K \lra \Om_C \lra \Om_{C/K} \lra 0
$$
(with $\Om_C$ projective).  Applying $\Hom(-,\O_C)$ yields the exact sequence
$$
\Hom(\Om_C,\O_C) \lra \Hom(f^*\Om_K, \O_C) \sr{\partial}{\lra} 
\Ext^1(\Om_{C/K},\O_C) \lra 0
$$
Given $x\in \Ext^1(\Om_{C/K},\O_C)$, choose $\beta \in \Hom(f^*\Om_K, \O_C)$
such that $\pd(\beta)=x$.  We then let $\zeta$, the pushout of $\eta$ along $\beta$,
be the extension corresponding to $x$.
$$
\xymatrix{
\eta: & 0 \ar[r] & f^*\Om_K \ar[r]\ar[d]^\beta 
& \Om_C \ar[r]\ar[d] & \Om_{C/K} \ar[r]\ar[d] & 0 \\
\zeta: & 0 \ar[r] & \O_C \ar[r] & \MG \ar[r] & \Om_{C/K} \ar[r] & 0
}
$$
Thus, by \cite{Ros}, the extension $\xi$ corresponds to the element
$\pd(\dltl) \in \Ext^1(\Om_{C/K},\O_C)$, where $\dltl$ is determined by the derivation 
on $K$.

There are natural isomorphisms $\Ext^1(\Om_{C/K},\O_C) \cong 
\Ext^1(\O_C, \Theta_{C/K}) \cong H^1(C,\Theta_{C/K})$, 
with $\Theta_{C/K} = \Om_{C/K}^\vee$, the tangent sheaf
of $C$ (the dual of $\Om_{C/K}$) (e.g., \cite{RH}, pages 234-5).
Thus, there is a natural map, which we also denote $\pd$,
$\pd:\Hom(f^*\Om_K,C) \lra H^1(C,\Theta_{C/K})$, 
which Buium~(\cite{Bui93}, p. 1395) calls the Kodaira-Spencer map of $C$.  
The element $\pd(\dltl) \in \Ext^1(\Om_{C/K},\O_C)$, or $H^1(C,\Theta_{C/K})$,
is the Kodaira-Spencer class of $C$.

From the sequence $\eta$, corresponding to the Kodaira-Spencer class of $C$,
one obtains a long exact sequence of cohomology groups.
$$
0 \lra \HH^0(C,\O_C) \lra \HH^0(C,\Omt_C) \lra \HH^0(C,\Om_{C/K})
\stackrel{\partial}{\lra}
\HH^1(C,\O_C) \lra \ldots
$$
Let us call $\dim_K \partial (\HH^0(C,\Om_{C/K}))$ the {\em Kodaira-Spencer} rank of
$C$, denoted $KS_C$.  Clearly, if the Kodaira-Spencer class of $C$ is trivial,
that is, if $\eta$ splits, then $KS_C = 0$, though the converse does not hold in 
general.  
For curves of genus 1, though, $\eta$ splits if and only if $KS_C = 0$.  
This follows, for example, 
from Atiyah's classification of rank 2 bundles over elliptic curves \cite{At}.
This yields a complete understanding of global $\tau$-forms on elliptic curves.

The following proposition is equivalent to Proposition~\ref{globdim} above.

\bpro
\label{cohom}
Let $C$ be a smooth projective curve.  
Then $1 \leq h^0(C,\Omt_C)\leq g + 1$.
Precisely, $h^0(C,\Omt_C) = g+1 - KS_C$.
In particular, if $C$ is defined over $k$, the field of constants, then 
$h^0(C,\Omt_C) = g+1 $.
\epro

\prf
The first two statements follow from remarks above.  The third follows from 
the fact that a curve defined over $k$ has trivial Kodaira-Spencer class.
(Since $H^1(C,\Theta_{C/K})$ classifies $TC$-torsors (e.g., 
see~\cite{JM}), this is equivalent to the fact that for such curves, 
$\tau C \cong_C TC$.)
\qed

\ignore{
To prove the third,
recall that the class $H^1(C,\Theta_{C/K})$ classifies $TC$-torsors (e.g., 
see~\cite{JM}), that and Buium~(\cite{Bui93}, p. 1396) proves that the 
{\em Kodaira-Spencer class}\
$\pd(\dltl)$ corresponds with the $TC$-torsor $\tau C$.  The torsor is  
trivial, that is, isomorphic to $TC$, if and only if $\pd(\dltl) =0$, which
holds exactly when the original sequence $\xi$ splits.  Putting everything together,
one gets that $\xi$ splits if and only if $\tau C \cong_C TC$, which happens exactly
when $C$ is defined over $k$, as desired.  This completes the proof.
\qed
}

\subsection{Definability results}

Recall that the {\em degree} of an invertible sheaf on a variety is defined to 
be the degree of the corresponding divisor.  We now observe that there is 
a (uniform) bound on the degree of a $\tau$-invertible sheaf on a curve.  

\bla
Let $V$ be a projective curve, $\MF$ a coherent sheaf on $V$.  Then for any 
invertible subsheaf $\MG$ of $V$, $\deg \MG \leq \hh^0(\MF) + g(C) - 1$.
\ela

\bpr
By the Riemann-Roch theorem, 
$\deg \MG = \hh^0(\MG) - \hh^0(\om \otimes \MG^{-1}) + g - 1 \leq 
\hh^0(\MG) + g - 1 $.  Since $\hh^0(\MG) \leq \hh^0(\MF)$, this yields the 
desired bound.
\epr

\bco
For any curve $C$, the degree of any \tis\ is $\leq 2g$.
\eco

The next result then follows by standard arguments.

\bco
For any curve $C$, the set of global $\tau$-forms on $C$ is definable,
as is the equivalence relation $\sim$ on these forms.  Further, this remains
true for families of curves (of fixed genus).  Finally, the dimension of 
the space of global $\tau$-forms, as a $K$-vector space, is uniformly definable.
\eco

Observe that, by Proposition~\ref{cohom}, this last claim is equivalent to the known fact
that the Kodaira-Spencer rank of a curve is definable.

\bibliography{jsl-smsets.bib}

\begin{thebibliography}{MMP06}

\bibitem[Ati57]{At}
M.~F. Atiyah.
\newblock Vector bundles over an elliptic curve.
\newblock {\em Proc. London Math. Soc. (3)}, 7:414--452, 1957.

\bibitem[Bui93]{Bui93}
Alexandru Buium.
\newblock Geometry of differential polynomial functions. {I}. {A}lgebraic
  groups.
\newblock {\em Amer. J. Math.}, 115(6):1385--1444, 1993.

\bibitem[GHL90]{GHL}
Sylvestre Gallot, Dominique Hulin, and Jacques Lafontaine.
\newblock {\em Riemannian geometry}.
\newblock Universitext. Springer-Verlag, Berlin, second edition, 1990.

\bibitem[Har77]{RH}
Robin Hartshorne.
\newblock {\em Algebraic geometry}.
\newblock Springer-Verlag, New York, 1977.
\newblock Graduate Texts in Mathematics, No. 52.

\bibitem[HI03]{HI}
E.~Hrushovski and M.~Itai.
\newblock On model complete differential fields.
\newblock {\em Trans. Amer. Math. Soc.}, 355(11):4267--4296 (electronic), 2003.

\bibitem[Hod93]{Hod}
Wilfrid Hodges.
\newblock {\em Model theory}, volume~42 of {\em Encyclopedia of Mathematics and
  its Applications}.
\newblock Cambridge University Press, Cambridge, 1993.

\bibitem[Hru95]{HJ}
Ehud Hrushovski.
\newblock {ODE}'s of order 1 and a generalization of a theorem of {J}ouanolou.
\newblock {\em Manuscript}, 1995.

\bibitem[Hru98]{HICM}
Ehud Hrushovski.
\newblock Geometric model theory.
\newblock In {\em Proceedings of the International Congress of Mathematicians,
  Vol. I (Berlin, 1998)}, number Extra Vol. I, pages 281--302 (electronic),
  1998.

\bibitem[HS94]{HS}
Ehud Hrushovski and Zeljko Sokolovi\'{c}.
\newblock Strongly minimal sets in differentially closed fields.
\newblock {\em Manuscript}, 1994.

\bibitem[HS99]{HrSc}
Ehud Hrushovski and Thomas Scanlon.
\newblock Lascar and {M}orley ranks differ in differentially closed fields.
\newblock {\em J. Symbolic Logic}, 64(3):1280--1284, 1999.

\bibitem[KN96]{KN96}
Shoshichi Kobayashi and Katsumi Nomizu.
\newblock {\em Foundations of differential geometry. {V}ol. {I}}.
\newblock Wiley Classics Library. John Wiley \& Sons Inc., New York, 1996.
\newblock Reprint of the 1963 original, A Wiley-Interscience Publication.

\bibitem[Mar00]{Mar03}
David Marker.
\newblock Manin kernels.
\newblock In {\em Connections between model theory and algebraic and analytic
  geometry}, volume~6 of {\em Quad. Mat.}, pages 1--21. Aracne, Rome, 2000.

\bibitem[Mat89]{Mat}
Hideyuki Matsumura.
\newblock {\em Commutative ring theory}, volume~8 of {\em Cambridge Studies in
  Advanced Mathematics}.
\newblock Cambridge University Press, Cambridge, second edition, 1989.
\newblock Translated from the Japanese by M. Reid.

\bibitem[Mil80]{JM}
James~S. Milne.
\newblock {\em \'{E}tale cohomology}, volume~33 of {\em Princeton Mathematical
  Series}.
\newblock Princeton University Press, Princeton, N.J., 1980.

\bibitem[MMP06]{MMP}
David Marker, Margit Messmer, and Anand Pillay.
\newblock {\em Model theory of fields}, volume~5 of {\em Lecture Notes in
  Logic}.
\newblock Association for Symbolic Logic, La Jolla, CA, second edition, 2006.

\bibitem[Pil02]{Pil02}
Anand Pillay.
\newblock Differential fields.
\newblock In {\em Lectures on algebraic model theory}, volume~15 of {\em Fields
  Inst. Monogr.}, pages 1--45. Amer. Math. Soc., Providence, RI, 2002.

\bibitem[Ros07]{Ros}
Eric Rosen.
\newblock Differentials over differential fields.
\newblock {\em Preprint}, 2007.
\newblock \texttt{arXiv:math.AC/0701508}.

\bibitem[Sha94]{Sha}
Igor~R. Shafarevich.
\newblock {\em Basic algebraic geometry. 1: Varieties in projective space}.
\newblock Springer-Verlag, Berlin, second edition, 1994.
\newblock Translated from the 1988 Russian edition and with notes by Miles
  Reid.

\bibitem[Wei94]{CW}
Charles~A. Weibel.
\newblock {\em An introduction to homological algebra}, volume~38 of {\em
  Cambridge Studies in Advanced Mathematics}.
\newblock Cambridge University Press, Cambridge, 1994.

\end{thebibliography}
\bibliographystyle{alpha}

\end{document}